\documentclass[12pt,a4paper]{amsart}
\usepackage{graphics}
\usepackage{epsfig}
\usepackage{graphicx}
\usepackage{enumerate, amssymb,stackrel,amsmath}
\usepackage[colorlinks]{hyperref}

\advance\hoffset 10mm \advance\textwidth40mm

\theoremstyle{plain}

\newtheorem{thm}{Theorem}[section]
\newtheorem{cor}[thm]{Corollary}
\newtheorem{lem}[thm]{Lemma}
\newtheorem{prop}[thm]{Proposition}
\theoremstyle{definition}
\newtheorem{defi}[thm]{Definition}
\newtheorem{defis}[thm]{Definitions}
\newtheorem{conj}[thm]{Conjecture}
\newtheorem{conv}[thm]{Convention}
\newtheorem{nota}[thm]{Setting}
\newtheorem{rem}[thm]{Remark}
\newtheorem{rems}[thm]{Remarks}
\newtheorem{exa}[thm]{Example}
\newtheorem{exas}[thm]{Examples}
\newtheorem{prob}[thm]{Problem}

\newtheorem{sit}[thm]{}
\newcommand{\rien}[1]{}

\newcommand{\Aut}{ \operatorname{{\rm Aut}}}

\def\ML{\mathop{\rm ML}}

\newcommand{\C}{\ensuremath{\mathbb{C}}}
\newcommand{\F}{\ensuremath{\mathbb{F}}}
\newcommand{\K}{\ensuremath{\mathbb{K}}}

\newcommand{\Q}{\ensuremath{\mathbb{Q}}}
\newcommand{\Z}{\ensuremath{\mathbb{Z}}}

\newcommand{\G}{\ensuremath{\mathbb{G}}}

\newcommand{\proj}{\ensuremath{\mathbb{P}}}

\newcommand{\kk}[1]{\bk^{[#1]}}

\newcommand{\bD}{{\bar D}}
\newcommand{\bX}{{\bar X}}
\newcommand{\bU}{{\bar U}}

\newcommand{\bS}{{\bar S}}
\newcommand{\hC}{{\hat C}}
\newcommand{\hF}{{\hat F}}

\newcommand{\hX}{{\hat X}}

\newcommand{\bC}{{\bar C}}
\newcommand{\bF}{{\bar F}}
\newcommand{\bY}{{\bar Y}}
\newcommand{\bZ}{{\bar Z}}

\newcommand{\tB}{{\tilde B}}
\newcommand{\tD}{{\tilde D}}

\newcommand{\tX}{{\tilde X}}

\newcommand{\hE}{{\hat E}}
\newcommand{\hD}{{\hat D}}

\newcommand{\cB}{{\ensuremath{\mathcal{B}}}}

\newcommand{\cF}{{\ensuremath{\mathcal{F}}}}
\newcommand{\cG}{{\ensuremath{\mathcal{G}}}}

\newcommand{\cO}{{\ensuremath{\mathcal{O}}}}

\newcommand{\cC}{{\ensuremath{\mathcal{C}}}}

\newcommand{\cX}{{\ensuremath{\mathcal{X}}}}

\newcommand{\p}{\partial}

\def\SAut{\mathop{\rm SAut}}
\def\kk{{\Bbbk}}

\def\kk{{\Bbbk}}

\def\PP{{\mathbb P}}

\def\pr{\mathop{\rm pr}}

\def\ML{\mathop{\rm ML}}

\newcommand{\AAutH}{ \operatorname{{\rm AAut}_{\rm hol}}}
\newcommand{\AutH}{ \operatorname{{\rm Aut}_{\rm hol}}}
\newcommand{\AutA}{ \operatorname{{\rm Aut}_{\rm alg}}}

\newcommand{\tC}{{\tilde C}}

\renewcommand{\epsilon}{\varepsilon}
\renewcommand{\phi}{\varphi}

\newcommand{\bnum}{\begin{enumerate}}
\newcommand{\enum}{\end{enumerate}}
\renewcommand{\emptyset}{\varnothing}
\newcommand{\brem}{\begin{rem}}
\newcommand{\brems}{\begin{rems}}
\newcommand{\erem}{\end{rem}}
\newcommand{\erems}{\end{rems}}
\newcommand{\bprob}{\begin{prob}}
\newcommand{\eprob}{\end{prob}}
\newcommand{\bexa}{\begin{exa}}
\newcommand{\bexas}{\begin{exas}}
\newcommand{\eexa}{\end{exa}}
\newcommand{\eexas}{\end{exas}}
\newcommand{\bdefi}{\begin{defi}}
\newcommand{\edefi}{\end{defi}}
\newcommand{\bdefis}{\begin{defis}}
\newcommand{\edefis}{\end{defis}}
\newcommand{\bcor}{\begin{cor}}
\newcommand{\ecor}{\end{cor}}
\newcommand{\blem}{\begin{lem}}
\newcommand{\elem}{\end{lem}}
\newcommand{\bconv}{\begin{conv}}
\newcommand{\econv}{\end{conv}}
\newcommand{\bconj}{\begin{conj}}
\newcommand{\econj}{\end{conj}}
\newcommand{\bprop}{\begin{prop}}
\newcommand{\eprop}{\end{prop}}
\newcommand{\bthm}{\begin{thm}}
\newcommand{\ethm}{\end{thm}}
\newcommand{\bnota}{\begin{nota}}
\newcommand{\enota}{\end{nota}}
\newcommand{\bsit}{\begin{sit}}
\newcommand{\esit}{\end{sit}}
\newcommand{\be}{\begin{eqnarray}}
\newcommand{\ee}{\end{eqnarray}}
\newcommand{\bproof}{\begin{proof}}
\newcommand{\eproof}{\end{proof}}
\def\ba{\begin{array}}
\def\ea{\end{array}}

\newcommand{\nlin}{\unitlength1mm\begin{picture}(0,9.25)
                      \put(0,0.75){\line(0,1){8.5}}
                     \end{picture}}
\newcommand{\slin}{\unitlength1mm\begin{picture}(0,0)
                      \put(0,-0.75){\line(0,-1){8.5}}
                     \end{picture}}

\newcommand{\vlin}[1]{\hspace{0.75mm}\unitlength1mm\begin{picture}(#1,0)
                      \put(0,0){\line(1,0){#1}}
                     \end{picture}\hspace{0.75mm}\rule[-3mm]{0mm}{4mm}}

\newcommand{\lin}{\vlin{8.5}}

 \newcommand{\mebox}{\unitlength1mm\begin{picture}(0,6)
   \put(-5,-2){\line(0,1){6}}
   \put(-5,-2){\line(1,0){6}}
   \put(1,4){\line(0,-1){6}}
   \put(1,4){\line(-1,0){6}}
   \end{picture}}

\newcommand{\cou}[2]{\unitlength1mm\begin{picture}(0,8)
   \put(0,0){\circle{1.5}}
   \put(0,3){\makebox(0,5)[b]{$#1$}}
   \put(0,-7){\makebox(0,4)[t]{$#2$}}
     \end{picture}
     \rule[-7mm]{0mm}{7mm}}

\newcommand{\crl}[2]{\unitlength1mm\begin{picture}(0,8)
   \put(0,0){\circle{1.5}}
   \put(-5,0){\makebox(0,5)[b]{$#1$}}
  \put(5,0){\makebox(0,5)[b]{$#2$}}
     \end{picture}
     \rule[-7mm]{0mm}{7mm}}

\newcommand{\cshiftup}[2]{\unitlength1mm\begin{picture}(0,9.25)
                      \put(0,10){\crl{#1}{#2}}
                     \end{picture}}
                     
\newcommand{\cshiftdown}[2]{\unitlength1mm\begin{picture}(0,9.25)
                      \put(0,-10){\crl{#1}{#2}}
                     \end{picture}}

\begin{document}
\title[Complete algebraic vector fields on affine surfaces]
{Complete algebraic vector fields on affine surfaces}

\author{Shulim Kaliman, Frank Kutzschebauch, and
Matthias  Leuenberger}
\address{Department of Mathematics,
University of Miami, Coral Gables, FL 33124, USA}
\email{kaliman@math.miami.edu}
\address{Department of Mathematics\\
University of Bern\\
Bern, Switzerland, }
\email{frank.kutzschebauch@math.unibe.ch}
\address{Department of Mathematics\\
University of Bern\\
Bern, Switzerland, }
\email{matthias.leuenberger@math.unibe.ch}

\begin{abstract}  Let $\AAutH (X)$ be the subgroup of the group $\AutH (X)$ of holomorphic automorphisms
of a normal affine algebraic surface $X$ generated by elements of flows associated with complete algebraic vector
fields.  Our main result is a classification of all  normal affine algebraic surfaces $X$ quasi-homogeneous under $\AAutH (X)$
in terms of the dual graphs of the boundaries $\bX \setminus X$ of their
SNC-completions $\bX$.   \end{abstract}

\maketitle

\date{}

\thanks{
{\renewcommand{\thefootnote}{} \footnotetext{ 2000
\textit{Mathematics Subject Classification:}
14R20,\,32M17.\mbox{\hspace{11pt}}
\\


This work was done during visits of the last two authors to the University of Miami and during a visit of the first author to the University of Bern.
We thank these institutions for their hospitality.The work of the last two authors was partially supported by Schweizerischer Nationalfond grant 200021-140235/1.}}}

{\footnotesize \tableofcontents}

\section{Introduction}  

\subsection{Motivation and general background}
 In the last decades  affine algebraic varieties and Stein manifolds with big (infinite-dimensional) automorphism groups have been
studied intensively. Several notions expressing the fact that the automorphisms group of a manifold is big have been proposed. 
Among the most important of them  are 
(algebraic) density property and holomorphic  flexibility with the former implying the latter.
Both density property and holomorphic  flexibility show that the manifold in question is an Oka-Forstneri\v c manifold. This important notion has also recently merged from  the intensive studies
around the homotopy principle which goes back to the 1930's and has had an enormous impact
on the development of Complex Analysis with a constantly growing number of applications (for definitions and more information we
refer the reader to \cite{For}). The newly emerged area of  Holomorphic Elliptic Geometry is devoted to the study of these properties.
 
In spite of the large number of examples of such highly symmetric  objects  their classification and the exact relations between all mentioned properties remain unclear
even in dimension 2. In particular, we do not know the description of Stein surfaces $X$ on which the group of holomorphic automorphisms $\AutH (X)$ acts transitively.
Since such transitivity is an automatic consequence of flexibility,  its study  is an important first step in our program of  {\bf Finding the exact relations between all known 
properties from Holomorphic  Elliptic Geometry for surfaces.}

In the algebraic case the analogous  question  of algebraic transitivity is classical, it was ``almost"  answered  in the papers of Gizatullin and Danilov \cite{Gi}, \cite{GiDa1}.  We need the following definition
to formulate their result.

\bdefi\label{quasi} We call a normal Stein (resp. affine algebraic) surface $X$  quasi-homogeneous with respect to a subgroup $G$ of the group of its holomorphic  (resp. algebraic)
automorphisms  if the natural action of $G$ has an open orbit
in $X$ whose complement is 
at most finite.
A normal affine algebraic surface is called quasi-homogeneous (without any reference to a group) if it is quasi-homogeneous with respect to 
the group $\AutA (X)$ of algebraic automorphisms\footnote{ When the complement to the open orbit is empty one has  transitivity.  However, there are
examples of smooth quasi-homogeneous surfaces for which the complements of the open orbits are not empty. In the case of surfaces over algebraically closed field of positive characteristic  they appeared already in the paper of
Gizatullin and Danilov \cite{GiDa} who also knew but did not publish such examples for characteristic zero. In a published form examples
of complex quasi-homogeneous surfaces with non-empty complements can be found  in a recent paper of Kovalenko \cite{Ko}.}.

\edefi

For convenience of the reader let us also recall the following

\bdefi
An SNC-completion (short for simple normal crossing)  of a normal affine algebraic surface $X$ is a normal complete algebraic surface $\bX$ such that
all  irreducible  components of the boundary curve 
$D= \bX \setminus X$ are smooth and all of its singularities are simple nodes (i.e., locally each of these singularities
is a transversal intersection of two smooth analytic branches). The boundary $D$ is contained in the smooth part of the surface $\bX$.
\edefi

Every quasi-homogeneous surface is either the two-dimensional torus $\C^* \times \C^*$, or  $\C \times \C^*$
or it   admits an SNC-completion
$\bX$ such that the dual graph $\Gamma$   of its boundary $\bX \setminus X$ (see Definition \ref{def.graph.1}) is a linear rational graph  \cite{Gi}, \cite{GiDa1} which can be always chosen
in the following standard form (so-called standard zigzag)
$$ \cou{C_0}{0}\lin\cou{C_1}{0}\lin\cou{C_2}{w_2}\lin\ldots\lin\cou{C_{n}}{w_{n}} \, \, \,  $$
where $n \geq 1$ and    $w_i \leq -2$ for $i=2, \ldots, n$.
A surface that admits such a completion will be called below a
Gizatullin surface\footnote{ In fact, for a normal affine surface to be a Gizatullin one it suffices to require that the dual graph of the boundary divisor is contractible
to a linear graph. 
The class of Gizatullin surfaces coincides also with the class of surfaces
with a trivial Makar-Limanov invariant, i.e., $\ML (X)=\C $ where $\ML (X)$ is the subring of the ring of regular functions on $X$ that consists
of all functions invariant under any (algebraic) $\G_a$-action on $X$}. Let $\SAut (X)$ be
the subgroup of $\AutA(X)$ generated by elements of $\G_a$-actions\footnote{Recall that
given a variety $V$ over a field $\kk$ an algebraic action of the additive group $\kk$ (resp. multiplicative group $\kk^*=\kk \setminus \{ 0 \}$)
on $V$ is called a $\G_a$-action (resp. $\G_m$-action.) In particular, in the case of $\kk=\C$ a $\G_a$-action is the same as an
algebraic action of the group $\C_+$ of complex numbers with addition as the group operation.)}  (also called algebraic $\C_+$ actions).
Then this subgroup possesses already an open orbit in $X$ whose complement is at most finite. Recall that  a $\G_a$-action can be viewed as the  flow
of a complete algebraic vector field on $X$ that is locally nilpotent.

More precisely,  a holomorphic vector field $\nu$ on a complex space $X$ is called complete if the solution of the ODE
$$\frac{d}{dt} \varphi (x,t) = \nu (\varphi (x,t)) \quad \varphi (x,0) = x$$

\noindent
is defined for all complex times $t \in \C$ and all initial values $x \in X$. The induced maps
$\Phi_t : X \to X$ given by $\Phi_t (x) = \varphi (x,t)$
yield the  flow of $\nu$ which is nothing but a one-parameter subgroup in the group $\AutH (X)$  of holomorphic
automorphism with parameter $t \in \C_+$ (equivalently, a holomorphic $\C_+$-action). 

One could wish to extend the quasi-homogeneity results to the analytic situation replacing locally nilpotent vector fields by complete holomorphic vector fields and $\SAut (X)$ by its holomorphic analogue,
but unfortunately the classification of complete holomorphic vector fields on Stein surfaces with sufficiently big automorphism groups
(even on $\C^2$)
seems still far out of reach. However, complete algebraic vector fields have been classified on $\C^* \times \C^*$ (actually on $(\C^*)^n$) by Anders\'en using Nevanlinna theory \cite{A'} and on $\C^2$ by Brunella  \cite{Bru04} using foliation theory.  
 It is important to remember  that when $X$ is an affine algebraic variety and $\nu$ is a complete algebraic vector field its  flow may well be non-algebraic. When
the  flow is algebraic then either $\nu$ is locally nilpotent and we have  a $\G_a$-action, i.e., an algebraic $\C_+$-action,  or $\nu$ is semi-simple and we have a 
$\G_m$-action, i.e., an algebraic $\C^*$-action.\footnote{It is worth reminding that an algebraic vector field $\nu$ on an affine algebraic variety $X$
is locally nilpotent if and only if for every element $f$ of the algebra $A$ of regular functions on $X$ there exists $n=n(f)$ such that $\nu^n (f)=0$. Similarly,
$\nu$ is semi-simple if and only if $A$ is a $\Z$-graded algebra $A=\bigoplus_{k\in \Z} A_k$ such that for every $f \in  A_k$ one has $\nu (f)=kf$. }

\bdefi\label{int.2} A holomorphic automorphism $\alpha$ of an algebraic variety $X$ will be called   {\it algebraically generated}  if $\alpha$ coincides with an element
$\Phi_t$ of the  flow of a complete algebraic vector field on $X$ as before.
The subgroup of the holomorphic automorphism group generated by such  algebraically generated automorphisms will be denoted by $\AAutH (X)$ 
and a normal affine algebraic surface $X$ quasi-homogeneous with respect to $\AAutH (X)$ will be called {\it generalized Gizatullin surface}. 
If a normal affine algebraic surface $Y$ admits two complete non-proportional algebraic vector fields  $\nu_1$ and $\nu_2$ 
(i.e., $f_1\nu_1 \ne f_2 \nu_2$ for any pair of nonzero regular functions $f_1$ and $f_2$ on $Y$)
then there is an open orbit of the natural $\AAutH (Y)$-action in $Y$. In what follows we call such $Y$ a  {\it surface with an open orbit }  (without mentioning $\AAutH (Y)$).
Of course, every  {\it generalized Gizatullin surface} 
is a surface with an open orbit.
\edefi

It is worth mentioning that the first examples of generalized Gizatullin surfaces
with discrete algebraic automorphism group $\AutA (X)$  were found by the
first and second author in \cite{KaKu1}, \cite{KaKu2}
 and they can be presented as hypersurfaces
 $\{ xp(x)+yq(y) +xyz=1 \}  \subset \C^3_{x,y,z}$ where the polynomials $1-xp(x)$ and $1-yq(y)$ have simple roots only  (none of these surfaces
 admits nontrivial algebraic $\G_a$ or $\G_m$-actions). Similar to the case of the torus for each such a surface $X$ an SNC-completion $\bX$ can be chosen
as a cycle. 

In this paper we deal with the following  first step of our general program. 

\subsection{Complete algebraic vector fields and  quasi-homogeneity by algebraically generated automorphisms}

The main result of our paper is the following classification of normal affine algebraic surfaces quasi-homogeneous under the group of algebraically generated automorphisms.

 \bthm\label{mthm}
 A  normal  complex affine algebraic surface $X$ is  generalized Gizatullin if and only if it admits an SNC-completion $\bX$  for which
 the boundary $\bX \setminus X$ is connected, consists of rational curves, and has a dual graph  $\Gamma$ that belongs to one of the following types\\
 \begin{enumerate}
 \item a standard zigzag or a linear chain of three $0$-vertices (i.e., Gizatullin surfaces and $\C \times \C^*$), \vspace{0.3cm}
\item circular graph with the following possibilities for weights\\[2ex] 
(2a) \hspace{1cm} $((0,0,w_1, \ldots , w_n))$ where $n \geq 0$ and $w_i \leq -2$, \\[2ex]
(2b) \hspace{1cm} $((0,0,w))$ with $-1\leq w \leq 0$ or $((0,0,0,w))$ with $w \leq 0$, \\[2ex]
(2c) \hspace{1cm} $((0,0,-1,-1))$;
\vspace{1cm} 
 \item $\hspace{1.5cm}\bigskip \cou{E \, \, \, \, \, \, \,  }{-1 \, \, \, \, \, \, \, \, \, }\nlin\cshiftup{\tilde
C_1}{-2}\slin\cshiftdown{\tilde C_2}{-2}\lin\cou{C_0}{w_0}\lin\cou{C_1}{w_1}\lin\ldots\lin\cou{C_{n}}{w_{n}} \hspace{0.5cm},\hspace{0.5cm} 
\begin{aligned} &
\end{aligned} $\vspace{1cm} 
\item $\hspace{1.5cm}
\bigskip
\cou{E \, \, \, \, \, \, \,  }{-1 \, \, \, \, \, \, \, \, \, }\nlin\cshiftup{\tilde C_1}{-2}\slin\cshiftdown{\tilde
C_2}{-2}\lin\cou{C_0}{-2} \hspace{4.5cm}\quad 
$\vspace{1cm} 
\item $\hspace{1.5cm}\bigskip
\cou{E \, \, \, \, \, \, \,  }{-1 \, \, \, \, \, \, \, \, \, }\nlin\cshiftup{\tilde C_1}{-2}\slin\cshiftdown{\tilde
C_2}{-2}\lin\cou{C_0}{w_0}\lin\cou{C_1}{w_1}\lin\ldots\lin\cou{C_{n}}{w_{n}}\lin\cou{ \, \, \, \, \, \, \, E'}{ \, \, \, \, \, \, \, k'}\nlin\cshiftup{\tilde
C_1'}{-2}\slin\cshiftdown{\tilde C_2'}{-2}
\quad,\hspace{1cm}
\begin{aligned} &
\end{aligned}$ \vspace{1cm} 
\item $\hspace{1.5cm}
\bigskip 
\cou{E \, \, \, \, \, \, \,  }{-1 \, \, \, \, \, \, \, \, \, }\nlin\cshiftup{\tilde C_1}{-2}\slin\cshiftdown{\tilde C_2}{-2}{\vlin{25}}    
\cou{ \, \, \, \, \, \, \, E'}{ \, \, \, \, \, \, \, k'}\nlin\cshiftup{\tilde C_1'}{-2}\slin\cshiftdown{\tilde C_2'}{-2} \hspace{1cm}, \hspace{1cm}
\text{for } k'\geq -1$ \vspace{1cm}

\end{enumerate}
\noindent  where $n\geq 0$, $w_i\leq -2$ for $i \geq 1$ and

{\rm ($\alpha$)}  in (3) either $n=0$ and $0\leq w_0\leq 2$ or $n \geq 1$ and $0\leq w_0\leq 1$;

{\rm ($\beta$) } in (5) either $n=w_0=0$ and  $k'\leq -1$ or $n \geq 1$, $0\leq w_0\leq 1$ and $k'\leq -3$ 
 or $n \geq 1$, $0\leq w_0\leq 1$, $w_n \leq -3$ and $k'=-2$;

{\rm ($\gamma$)} in (2a) $\Gamma$ is a subgraph of a graph contractible to the cycle  $((0,0,0,0))$.

 \ethm
 
  It should be mentioned that the second author announced a formulation of Theorem \ref{mthm} in the survey \cite{Ka19}
but unfortunately in that formulation conditions ($\alpha$) and ($\beta$) were stated incorrectly which results in
some extra graphs that do not appear as the dual graphs of SNC-curves in smooth complete
surfaces.

In the framework of our general program Theorem \ref{mthm} provides us with a list\footnote{It should be emphasized that in general two surfaces
with the same boundary graph as one of those which appear  in Theorem \ref{mthm} are not necessarily homeomorphic. Even in the case of the same topology such surfaces may admit families with non-isomorphic members (we
do not know if a similar fact holds in the analytic setting).   Furthermore, homogeneity of such surfaces with respect to the $\AAutH$-action does not guarantee algebraic density property.
For instance, the hypersurface $\{ xp(x)+yq(y) +xyz=1 \}  \subset \C^3_{x,y,z}$ mentioned before does not have it, though it has another
positive feature -  the algebraic volume density property \cite{KaKu2}.}
of affine algebraic surfaces on which one can look for  holomorphic flexibility, algebraic
density property, or even classification of  all complete algebraic vector fields. A  case of particular surfaces has been already worked out by the third author in \cite{L}. 
Following our general program the algebraic density property for certain classes of  Gizatullin surfaces has been recently established  by Andrist, Poloni and the second author in \cite{AKP} and  by Andrist in \cite{A}. 

The second important result of our paper is one of the essential ingredients in the proof of our  main result.

\bthm \label{mthm1i} Let $X$ be a normal affine algebraic surface which admits a nonzero complete algebraic vector field.  Then  either:

(1) all complete algebraic fields share the same rational first integral (i.e., there is a rational map $f: X \dashrightarrow B$
 onto a smooth curve $B$ such that all complete algebraic vector fields on $X$ are tangent to the fibers of $f$), or

(2) $X$ is a rational surface with an open orbit and, furthermore, for every complete algebraic vector field $\nu$ on $X$ there is  a regular function $f: X \rightarrow \C$ 
(depending on $\nu$) with general fibers isomorphic to $\C$ or $\C^*$ such that
the flow of $\nu$ sends fibers of $f$ to fibers of $f$ and  $f(X)$ is either $\C$ or $\C^*$.
\ethm

The fact that the flow sends fibers to fibers can be reformulated as follows: there is a complete vector field $\nu_0$ on $\C$ such that $f$ is $\C_+$-equivariant with
respect to the  flows of $\nu$ (resp. $\nu_0$) acting on $X$ (resp. $\C$). When $\nu_0$ is trivial then $f$ is again a rational
first integral of $\nu$ (which is regular in this case)

In the special case of  $X = \C^2$ Theorem  \ref{mthm1i} was proven by Brunella \cite{Bru04} and our original proof followed his approach
(see Remark \ref{}).
 However,  after finishing the first version of the manuscript we were informed about the paper of Guillot and Rebelo \cite{GR}.
They proved a very general analogue of Theorem  \ref{mthm1i}  for  so called semicomplete meromorphic  vector fields 
on complex surfaces in which the conclusion about
 the isomorphism type of a general fiber of $f$ in (2) is replaced by the following: the completion of such a fiber is either a rational curve
or an elliptic one. Complete algebraic fields on affine algebraic surfaces are, of course, semicomplete, and
it is easier to extract Theorem  \ref{mthm1i} from the Guillot-Rebelo result which we do  in the present version.

The paper is organized as follows.  In Section \ref{theorem1.3} we start the proof of a crucial part, the case of absence of  a rational first integral,  of Theorem   \ref{mthm1i} based on the Guillot-Rebelo theorem.

In Section \ref{dualgraphs} we remind  some  facts from \cite{FKZ} (see also \cite{Dai1}, \cite{Dai2}) about weighted
dual graphs of algebraic curves contained in
smooth surfaces while in
Section \ref{P1fibrations} we present some basic results on $\PP^1$-fibrations on $\bar X$ extending $\C$- or $\C^*$-fibrations on $X$. 

In Section \ref{non-affine}, with these preliminaries in hand we can finish the case of absence of  a rational first integral  from  Theorem   \ref{mthm1i}. 

Section \ref{rational} is devoted to the geometrical description of rational
first integrals of complete algebraic vector fields on smooth semi-affine surfaces (whenever such integrals exist) which together with  the results from  Section \ref{theorem1.3} and Section \ref{non-affine}
allows us to complete the proof of Theorem \ref{mthm1i}.

In Section \ref{existinggraphs} we show that the graphs  listed in Theorem \ref{mthm} appear as the dual graphs
 of SNC-curves with rational irreducible components contained in smooth rational surfaces.

In Section \ref{boundary} we show there is no graph different from those described in  (1)-(6) that can serve as the dual graph of
a boundary of a surfaces quasi-homogeneous under $\AAutH (X)$.

In Section \ref{proof} we establish that any surface with a boundary graph as in (1)-(6) is indeed quasi-homogeneous under $\AAutH (X)$
which  concludes  the proof of  Theorem \ref{mthm}.

In the last Section  we show that  some of  the surfaces in Theorem \ref{mthm} are in fact $\AAutH (X)$-homogeneous.
In particular, this class includes all surfaces of type (2) and all known Gizatullin surfaces that are not homogeneous with respect to the natural $\Aut (X)$-action described in  \cite{Ko}.


\section{Existence of Riccati fibration}\label{theorem1.3}

\bdefi\label{GR.semi}  
A semi-affine surface is an algebraic surface $X$ that admits a proper birational morphism onto an affine algebraic surface $S$
(which in complex analysis is nothing but the Remmert reduction). 
For instance, resolution of singularities of
a normal affine surface $S$ leads, of course, to a smooth semi-affine surface $X$. 

\edefi

\brem\label{GR.boundary}   Note that  the boundary divisor $D=\bX \setminus X$ for some
SNC-completion $\bX$ of $X$ can be also viewed as the boundary divisor of some SNC-completion $\bS$ of $S$
and the intersection matrix of the irreducible components of $D$ in $\bX$ is the same as the similar matrix for $D$ in $\bS$.
\erem

 We shall use the following standard fact.

\bprop\label{standard} 
Algebraic vector fields on $S$ lift to algebraic vector fields on $X$ and complete ones lift to complete ones. Vice versa, complete vector fields on $X$ correspond to complete vector fields on $S$. 
\eprop

\bproof For simplicity let us deal with complete algebraic fields only. Consider the  flow of such a field. The singularities of $S$ are fixed by this flow.
Hence, it induces a 
holomorphic $\C_+$-action on the tangent cones at these points, and, therefore, on the result  $\tilde S$ of blowing $S$ up at these singularities (since the
exceptional curves in $\tilde S$ are projectivizations of tangent cones).
Then the action on $\tilde S$ induces a holomorphic $\C_+$-action on its normalization. Since a minimal resolution of surface singularities involves
only a sequence of blowing-up and normalizations \cite{Li} we get a holomorphic $\C_+$-action on $X$ and, thus, the desired complete algebraic vector field on $X$.

Vice versa, given a complete algebraic vector field on $X$ we note that its  flow preserves complete curves contained in $X$
since the number of such curves in a semi-affine surface is finite. Thus, one can push this flow down to the regular part of $S$
and extend the holomorphic $\C_+$-action to $S$ by Hartogs' theorem. This yields a complete algebraic vector field on $S$.
\eproof

Therefore, to consider complete algebraic vector fields on normal affine
surfaces is the same as to consider such fields on smooth semi-affine surfaces. 

For convenience of the reader we recall  the following

\bdefi\label{5.10}
A foliation $\cF$ on a smooth complex surface $\bX$ is given by an open
covering $\{ U_j \}$ of $\bX$ and holomorphic vector fields $\nu_j \in H^0(U_j, T\bX )$
with isolated zeros such that
$$\nu_i = g_{ij} \nu_j \, \, \, {\rm on} \, \, \, U_i \cap U_j$$ for invertible
holomorphic functions $g_{ij} \in H^0 (U_i \cap U_j, \cO^*_{\bar X} )$ where $\cO_{\bar X}^*$
is the sheaf of invertible functions. Gluing orbits
of $\{ \nu_j \}$ one gets leaves of the foliation $\cF$. The singular set
${\rm Sing} \, (\cF )$ is the discrete subset of $\bX$ whose intersection with
each $U_j$ coincides with the zeros of $\nu_j$. The cocycle $\{ g_{ij} \}$ defines
a holomorphic line bundle $K_{\cF}$ which is called the  canonical bundle
of the foliation $\cF$. 
\edefi

One of the main facts about foliations which is of use for us is the following consequence of the Camacho-Sad formula (e.g., see \cite{Bru00}).

\bprop\label{CS} Let $\cF$ be a foliation on a compact smooth surface $\bX$ and $C$ be an $\cF$-invariant closed curve
that does not contain singularities of $\cF$. Then the self-intersection $C \cdot C$ of $C$ is zero.

\eprop

\begin{nota}\label{notaresol}
Let $X$ be a smooth semi-affine smooth  algebraic surface
equipped with a nonzero complete algebraic vector field $\nu$
and an SNC-completion $\bX$ (i.e.,  $\bD = \bX \setminus X$ is an SNC-divisor in the smooth projective surface $\bX$). The algebraic vector field $\nu$ generates a
rational vector field $\bar \nu$ on $\bX$. Note also that $\nu$ generates a foliation $\cF$ on $X$ which extends
to the foliation $\bar \cF$ on $\bX$ generated by $\bar \nu$.  

\end{nota}

As we mentioned before, complete algebraic vector fields on $X$ form a subset in the wider class of semi-complete vector fields. Therefore,  in our case Theorem B in \cite{GR} 
implies the following.

\bthm\label{GR} In Setting  \ref{notaresol}  either

{\rm (1)} $\bar \nu$ has a first integral \footnote{I.e., a surjective rational map $\bX - \to C$ into a curve $C$
whose fibers are tangent to $\bar \nu$.}, or

{\rm (2)} up to a birational transformation of $\bX$ the field $\bar \nu$ is holomorphic, or

{\rm (3)} there is a morphism $f: \bX \to B$ into a complete rational or elliptic curve $B$ with  rational or elliptic general fibers 
such that for some vector field $\mu$ on $B$ one has $f_* (\bar \nu ) = \mu$ (the  flow of $\bar \nu$ preserves the fibration $f$).

\ethm

\bdefi\label{RT} A morphism $f : \bX \to B$ as above with rational (resp. elliptic) general fibers is called a Riccati (resp. Turbulent) fibration adapted to  the foliation $\bar \cF$.
\edefi

\bprop\label{GR.step}  In Setting \ref{notaresol} let $X$ be a semi-affine smooth surface and $f: \bX \to B$ be a fibration adapted to $\bar \cF$ as in
Theorem \ref{GR} (3).
Suppose that $\bar \nu$ does not have a rational first integral. Then $f$ is a Riccati fibration and
its restriction to $X$ yields a morphism
with general fibers isomorphic either to $\C^*$ or to $\C$. Furthermore, $B$ is a rational curve. 

\eprop

\bproof  Since $X$ is a semi-affine surface it may contain only a finite
number of complete curves. Hence, a general fiber $F=f^{-1}(b_0)$ of $f$ (which is isomorphic by definition either  to $\PP^1$ or to a complete elliptic curve) must meet the union $C$ of the components of boundary divisor
$D = \bX \setminus X$ that are not contained in the fibers of $f$. The absence of first integrals implies that $\bar \nu$ is not tangent to $F$ but it is tangent to $C$ since $\nu$ is a complete field.
In particular, $C$ consists of leaves of the foliation $\bar \cF$.

Furthermore, $\mu$ from Theorem \ref{GR} (3) (which is complete by construction) must be nonzero. Hence, removing from $B$ a finite number of points corresponding to non-general values of $f$ we get a curve $B^*$
which is an orbit of the  flow of $\mu$. Thus, if $B$ is rational then $B^*$ is isomorphic either to $\C$ or to $\C^*$ while in the case of an elliptic curve
$B^*=B$ since the lift of a vector field $\mu$  (which is not identically zero) 
to the universal cover $\C$ of $B$ may be only a nonzero constant field, i.e., $\mu$ does not vanish at any point of $B$. 

Note that foliation $\bar \cF$ generates a representation of the fundamental group $\pi_1 (B^*, b_0)$ into the group $\Aut_{\rm hol} (F)$ of holomorphic automorphisms of $F$ which must
preserve the finite nonempty subset $F \cap C$ (because the leaves of the foliation are invariant under monodromy).
 
To get a contradiction assume that  the curve $F \setminus C$ is of general type
(i.e., either $F$ is elliptic or $F \simeq \PP^1$ and $F \cap C$ consists at least of three points).
Then the image of $\pi_1(B^*, b_0)$ under the representation is contained in a finite subgroup $G \simeq \Aut (F \setminus C)$ of $\Aut_{\rm hol} (F)$.
Consider the unramified cover $\tB^* \to B^*$ 
for which $G$ plays the role of the Galois group,
the fibered product $\tX^* = X \times_B \tB^*$, the preimage $\tC^*$ of $C$ in $\tX^*$ and the natural projection $\tilde f : \tX^* \to \tB^*$.
Note that every irreducible component of $\tC^*$ is now a section of $\tilde f$ and  the vector field $\tilde \nu$ 
on $\tX^*$ induced by $\nu$ is tangent to $\tC^*$ and sends every fiber of $\tilde f$ to a fiber of $\tilde f$.

{\em Claim}. For the first statement it suffices to show that there is a Zariski dense open subset $U$ of $\tB^*$ and  $F\setminus C \subset F^*\subset F$ such that
(a) $F^*$ is a curve of general type and (b) $\tilde f^{-1} (U)$ is naturally isomorphic to $U\times F$ over $U$.

Indeed,  then $\tilde \nu$ is tangent to the curve $U\times (F \setminus F^*)\subset \tC^*$ (i.e.,
 the flow of $\tilde \nu$ transforms $U\times F^*$ into itself) and since $\tilde f _*(\tilde \nu)$ is a vector field on $U$ 
and $F^*$ is of general type for every element $U\times F^*\to U\times F^*$ of the flow the second coordinate
function must yield the identity map $F^* \to F^*$. That is, the natural projection $\tau : \tilde f^{-1} (U)\to F$ is a rational first integral for $\tilde \nu$.
Consider  now a rational function $h$ on
$F$ that has pole at every point of $F\cap C$. Let $h_1, \ldots , h_m$ be the orbit of $h$ under the natural action
of $G$. Since $C\cap F$ is invariant under the action of $G$ we see that each $h_i$ has a pole at every point of $F\cap C$.
In particular, the product $h_1\cdot \ldots \cdot h_m$ is a non constant  $G$-invariant
 meromorphic function on $F$. Taking the composition of this product with $\tau$
one gets a rational function $\tilde h$  on $\tX^*$ invariant with respect to the natural $G$-action. Hence, $\tilde h$ induces a rational function on $\bX$
which is a well-defined  first integral of $\bar \nu$ contrary to the assumption. This concludes the proof of the Claim.

Let us now check conditions (a) and (b).
If $F=\PP^1$ then $\tX^*$ is a trivial $\PP^1$-bundle over $\tB^*$ since $\tB^*$ is affine, i.e., we have (b) with $U=\tB^*$.
Note that $\tC^*$ consists at least of three irreducible components.
Choosing the isomorphism $\tX^* \simeq \tB^*\times \PP^1$ appropriately one can suppose that
these three sections can be viewed as the curves $\tB^*\times \infty$, $\tB^*\times 0$ and $\tB^*\times 1$ in $\tB^* \times \PP^1$.
Thus, letting $F^*=\PP^1 \setminus \{ 0,1,\infty\}$ we get (a) which concludes the proof  of the first statement in the case of $F=\PP^1$.

Now let us deal with the case when $F$ is elliptic. Consider the generic fiber of $\tilde f$. It is an elliptic curve $E$ over the field $\K$ of rational functions on $\tB^*$
and (being sections of $\tilde f$) the irreducible components of $\tC^*$ can be treated as closed points in $E$.
Recall that $E$ can be viewed as a smooth cubic in $\K\PP^2$ given by an equation in the Weierstrass form 
$\breve y^2=\breve x^3+ a_4\breve x \breve z^2 +a_6 \breve z^3$ 
where $a_4, a_6 \in \K$ and $(\breve x : \breve y : \breve z)$ is a homogeneous coordinate system on $\K\PP^2$ (e.g., see \cite{SchSh}).
In particular, $a_4=a_4(\tilde b)$ and $a_6=a_6 (\tilde b)$ are rational functions on $\tB^*$. Since the flow of $\tilde \nu$ transforms
the fibers of $\tilde f$ into the fibers of $\tilde f$ we see that for general $\tilde b_0$ and $\tilde b$ in $\tB^*$ the elliptic curves in $\C\PP^2$
given by the equations $\breve y^2=\breve x^3+ a_4(\tilde b_0)\breve x \breve z^2 +a_6(\tilde b_0) \breve z^3$ and 
$\breve y^2=\breve x^3+ a_4(\tilde b)\breve x \breve z^2 +a_6(\tilde b) \breve z^3$ are isomorphic.
Hence they have the same $j$-invariant (e.g., see \cite{SchSh}) which implies that, fixing $\tilde b_0$ one has $a_4(\tilde b)=a_4 (\tilde b_0)u(\tilde b)^4$
and $a_6(\tilde b)=a_6 (\tilde b_0)u(\tilde b)^6$ where $u(\tilde b)$ depends on $\tilde b$.  Assuming for simplicity that
$a_4 \ne 0 \ne a_6$ we observe that $u^2(\tilde b)$ is a rational function on $\tB^*$
and, thus, dividing by its powers we can suppose that for every $\tilde b$ in a nonempty Zariski open subset $U$ of $\tB^*$ the fiber $\tilde f^{-1} (\tilde b)$
is given by an equation $\breve y^2=\breve x^3+ a_4\breve x \breve z^2 +a_6 \breve z^3$  where $a_4$ and $a_6$ are now complex numbers
(i.e., we have (b)).
Since $E$ acts transitively on itself we can suppose that one of the irreducible components of $\tC^*$ corresponds to the point $(0:0:1)$
and that every fiber of $\tilde f^{-1} (\tilde b) \setminus \tC^*, \, \tilde b \in U$ is a curve contained in the affine curve $F^* \subset \C^2$
given by $ y^2= x^3+ a_4 x  +a_6 $ which yields (a). Thus, we have the first statement.

If we assume that $B$ is elliptic then it was shown already that $\mu$ and, therefore, $\bar \nu$ do not vanish. 
That is, every irreducible component of $C$ has no singularities of $\bar \cF$.  
This implies that these components are actually connected components of $C$ and furthermore $C$ coincides with the boundary divisor $D$ (indeed, otherwise
$C$ contains singularities of $D$ which must be automatically singularities of $\bar \cF$).
By Proposition \ref{CS} $D \cdot D=0$.  However,  by Remark \ref{GR.boundary} $D$ can serve also as a boundary of an SNC-completion of an affine surface.
In particular, it is connected (by the Lefschetz hyperplane section theorem) and being a support of an ample divisor (Theorem  2 in \cite{Good}) it has a positive self-intersection number $D \cdot D >0$. 
This contradiction implies that $B$ cannot be elliptic and, thus, $B$ is rational. 
\eproof

\brem This argument showing that $F\setminus C$ cannot be hyperbolic essentially appeared already  in the proof of  \cite[Lemma 1]{Bru04}
in which Brunella showed that the foliation induced by a complete vector field
on $\C^2$ cannot be a Turbulent one.
\erem

\bcor Under the assumption of Proposition \ref{GR.step} the surface $X$ is also rational.
\ecor

\bproof This follows from the fact that $B$ and general fibers of the morphism $f : \bX \to B$ are rational.
\eproof

\bprop\label{GR.other.step} The case (2) of Theorem \ref{GR} occurs only if one of cases (1) or (3) occurs as well.

\eprop

\bproof Suppose that $\psi : \bX - \to Y$ is a birational map of $\bX$ into a smooth projective surface $Y$ which transforms $\bar \nu$
into a holomorphic field $\mu$ on $Y$.

Note that the image $T$ of $D$ under $\psi$ cannot be a point. Indeed, by \cite{Har} both $\bX$ and $Y$ can be obtained via contractions $\alpha : Z \to \bX$ and $\beta : Z \to Y$ from
the same smooth surface $Z$. Note that for $E=\alpha^{-1}(D)$ the surface $Z \setminus E$ is semi-affine since $\bX \setminus D$ is. In particular, one has
a birational morphism $\gamma : Z \setminus E \to S$ into an affine surface $S$. By construction the rational map $\gamma \circ \beta^{-1} : Y \setminus T -\to S$ is birational
and, furthermore, it is regular by Hartogs' theorem. Thus, $Y \setminus T$ is semi-affine and $T$ is a divisor.

{\em Claim.} The Kodaira dimension of $Y$ cannot be non negative. 

Assume the contrary. Then by \cite[Ch. 6, Prop. 5]{Bru00} $\mu$ does not vanish and either it generates a foliation on $Y$ with elliptic curves
as general leaves, or $Y$ is a torus and $\mu$ generates a Kronecker foliation (for definition see \cite{Bru00}) on it. The former case is impossible. Indeed,  a general leave $G$ of the foliation $\cG$ induced by $\mu$
contains an isomorphic image of an integral curve of $\nu$ under $\psi$ (because the rational map $\psi$ is regular outside codimension 2).  Such an integral curve
is automatically isomorphic to $\C$ or $\C^*$ since $\nu$ is a complete vector field on the semi-affine surface $X$. Thus, $G$ is rational and we need to
consider the case of the Kronecker foliation only. 

By definition any Kronecker foliation $\cG$ has no algebraic curves among its 
leaves. Since $\mu$ is nowhere vanishing but every complete vector field on $\PP^1$ has zeros we see that
a general leaf $G$ of $\cG$ is isomorphic to either $\C$ or $\C^*$. If $G \simeq \C^*$ then $G$ is an image of a general integral curve of $\nu$.
Since these integral curves do not meet $D$ we conclude that $G$ does not meet $T$. This implies that  the algebraic curve $T$ is tangent to the 
non-vanishing field $\mu$ and, thus, it consists
of leaves of $\cG$. A contradiction. 

Thus, $G \simeq \C$. Since $T$ is not a leaf it must meet $G$ and furthermore $G\cap T$ is a singleton (indeed otherwise, a general integral curve of $\nu$, which is naturally isomorphic
to $G \setminus T$ is hyperbolic).  In fact compactness of $T$ implies that it meets every leaf of $\cG$ at one point. 
Applying the  flow of  the nowhere-vanishing $\mu$ to $T$ we see that $Y$ is biholomorphic to $T\times \C$ contrary to the fact that the torus $Y$ is a compact surface.
This concludes the proof of the Claim.

Thus, $Y$ is of negative Kodaira dimension.  Then according to \cite[Ch. 6, Prop.6]{Bru00} either $\mu$ is tangent to the fibers of the Albanese fibration
(and, thus, has a first integral), or $Y$ can be chosen as a quadric  and $\mu = \mu_1 \oplus \mu_2$ where $\mu_i$ is a field on a factor
in the direct product  $Y=\PP^1 \times \PP^1$ (and we have a natural projection $Y \to \PP^1$ as an adapted Riccati fibration), or it is the
suspension of a representation $\rho : \pi_1( {\rm Alb}(Y)) \to \Aut (\PP^1 )$ (and again the Albanese morphism $Y \to {\rm Alb}(Y)$ is  an adapted Riccati fibration).
This yields the desired conclusion. 
\eproof

Combining Propositions \ref{GR.step} and \ref{GR.other.step}  with Proposition \ref{standard} we get the following.

\bthm\label{GR.theorem} Let $\sigma$ be a  complete algebraic vector field on a normal affine surface $Z$ such that $\sigma$ does not
admit a rational first integral. Then there is a morphism $g: Z \to B$ onto a curve  $B$ isomorphic either to $\C^*, \C$ or to $\PP^1$ whose general fibers are isomorphic to $\C$ or $\C^*$ so that
the  flow of $\sigma$ sends fibers of $g$ to fibers of $g$. 
\ethm

In fact we will show in Section \ref{non-affine} (Theorem \ref{KKtheorem})  that in case that  the morphism $g$ is surjective to $\PP^1$, there is another morphism $g_1: Z \to B$ with affine base ($B$ is isomorphic to $\C$ or $\C^*$) whose general fibers are isomorphic to $\C$ or $\C^*$ so that
the  flow of $\sigma$ sends fibers of $g_1$ to fibers of $g_1$.

\brem\label{}
Let us briefly discuss the alternative proof in our original preprint based on the fundamental results of foliation theory for surfaces due
to Suzuki \cite{Su77a}, \cite{Su77b}, McQuillan \cite{McQ}, \cite{McQ08}, and Brunella \cite{Bru04}. 
Suppose that we are in Setting \ref{notaresol}. 
 By a theorem of Seidenberg \cite{Sei} after additional blowing-up $\hX \to \bX$ one can suppose that the induced foliation $\hat \cF$
on $\hX$ has reduced singularities only (see Definition \ref{new.d1} below).

Consider  the canonical bundle $K_{\hat \cF}\in {\rm Pic} (\hX) \otimes \Q$ of $\hat \cF$ as in Definition \ref{5.10}.
One can define the Kodaira dimension of $\hat \cF$ as the  Kodaira-Iitaka dimension of
$K_{\hat \cF}$. \footnote{That is,

$${\rm kod} (\hat \cF ) = \lim{\rm sup}_{n \to + \infty} {\frac{\log \dim H^0 (\bX, K_{\hat \cF}^{\otimes n})}{\log n }}.$$}
In the case when there is no  rational first integral  for our complete field $\nu$ the crucial result of
McQuillan implies that the Kodaira dimension $kod (\hat \cF)$ is either 1 or 0. For $kod (\hat \cF)=1$
a theorem of McQuillan yields an adapted Riccati fibration for $\hat \cF$. Brunella showed that this fibration induces a $\C$- or $\C^*$-fibration on $X$
whose fibers are transferred by the flow of $\nu$ into similar fibers (his proof is written for $\C^2$ but it works in general case as well).

For $kod (\hat \cF)=0$ one has so-called McQuillan contraction $\theta : \hX \to \hX'$ into a normal projective surface $\hX'$ so that $\hat \cF$ induces a foliation
$\hat \cF'$ on $\hX'$. Furthermore, there is a finite morphism $\rho : \hX' \to Y$ into a smooth surface $Y$ ramified at singular points of $\hX'$ only and such that
it transforms  $\hat \cF'$ into a foliation $\cG$ on $Y$ generated by a global holomorphic vector field $\mu$ on $Y$. We managed to show that
$\hX'$ is in fact smooth and, therefore, $\rho$ is an isomorphism. Then arguing as in Proposition \ref{GR.other.step} we showed that in the
absence of first integrals there must be a Riccati fibration on $\cX'=Y$ with desired properties.

\erem


\section{Preliminaries about dual graphs}\label{dualgraphs}

In this section we discuss some facts about weighted dual graphs of algebraic curves contained
in smooth complete surfaces (e.g., see \cite{FKZ}).

\bdefi\label{def.graph.1}
Let $D$ be a closed curve contained in a smooth complete algebraic surface $\bX$ such that all of its singularities are simple nodes (i.e., locally each of these singularities
is a transversal intersection of two smooth analytic branches). Then we can assign
the so-called weighted dual graph $\Gamma$ such that 

(1) its vertices are in bijective correspondence with the irreducible components of $D$;

(2) every singularity of $D$ corresponds to an edge that joins vertices corresponding to irreducible components $C_1$ and $C_2$
that contain this singularity (if $C_1=C_2$ then this edge is a loop);

(3) each vertex is equipped with a weight that is the integer equal to the self-intersection number $C^2$ of the corresponding component
$C$ in $\bX$.

(4) We also say that $D$ is of simple normal crossing type if all components are smooth. This implies that the dual graph $\Gamma$ does not have loops.
\edefi

\bconv\label{conv.graph.1}
From now on we identify the vertices of $\Gamma$ and the corresponding components of $D$ and denote them by the same letters. 
Furthermore, we may treat a curve contained in $D$ as a subgraph of $\Gamma$ and vice versa.

\econv

Recall that the valency of a vertex $C \in \Gamma$ is the number of edges adjacent to this vertex (with each loop counting twice) and the vertices 
joined with $C$ by edges are called the neighbors of $C$.
The vertex is an end vertex
(resp. linear vertex, resp. branch point) if its valency is 1 (resp. 2; resp. at least 3). The graph is called linear or a chain if it does not contain
branch points but contains an end vertex. We use notation as $C_1+C_2 + \ldots +C_n$ to denote such a chain with $n$ vertices
in the natural order. If the weight of each $C_i$ is $w_i$  we shall also use notation $[[w_1, w_2, \ldots , w_n]]$ for this chain.  A graph without branch
points and end vertices is called circular. In this case we write $((w_1, w_2, \ldots , w_n))$ for the weights of this graph in a counterclockwise order. 
For any subgraph $\Gamma_0$ of $\Gamma$ the notation $\Gamma \ominus \Gamma_0$ will be used to denote the graph obtained from $\Gamma$ by removing 
each vertex $C \in \Gamma_0$ and its adjacent edges. 

If $C$ is a rational irreducible component of $D$ with self-intersection $k$ we call $C$ a  $k$-vertex. If $C$ is a $(-1)$-vertex 
with valency at most 2 in $\Gamma$ then it can be contracted
and the image $D'$ of $D$ is still a curve with nodes as singularities (unless $C$ with a loop is a connected component of $\Gamma$ -
a case which we do not consider).
The graph $\Gamma'$ of $D'$  in the smooth resulting surface can be obtained from $\Gamma \ominus C$ by
joining the distinct neighbors of $C$ by an edge and increasing their weights by 1 (we call such replacement of $\Gamma$
by $\Gamma'$ a blowing down). 

A graph $\Gamma$ is contractible if  it  can be reduced to an empty graph by
a sequence of blowing downs (i.e $D$ can be blown down to a 
smooth point of a resulting surface).  We call $\Gamma$ minimal if it does not contain $(-1)$-vertices
different from branch points.

Let $z \in D$, $\sigma : \hX \to \bX$ be the monoidal transformation of $\bX$ at $z$, and $D''=\sigma^{-1}(D)$. Then the form of the
dual graph $\Gamma''$ of $D''$ in $\hX$ depends on whether $z$ is (a) a smooth point of $D$ (and in particular
$z$ is a smooth point of an irreducible component $C$ of $D$) or (b) a double point of $D$,
i.e., $z\in C_1 \cap C_2$ where $C_1$ and $C_2$ are irreducible components of $D$. In case (a) $\Gamma''$ is obtained from $\Gamma$
by creating a new vertex of weight $-1$, joining it with $C$, and reducing the weight of $C$ by 1. This change of $\Gamma$
will be called an outer blowing up. In case (b) $\Gamma''$ is obtained from $\Gamma$ by removing the edge between
$C_1$ and $C_2$, reducing their weights by 1, and joining them by edges with a new vertex of weight $-1$.
Such a change of $\Gamma$ will be called an inner blowing up.

If  a graph $\Gamma_2$ can be
obtained from a graph $\Gamma_1$ by a sequence of blowing up and blowing down then we call this procedure a
reconstruction of $\Gamma_1$ into $\Gamma_2$. Let us give some examples of reconstructions. If $C$ is of non negative weight in a graph $\Gamma \ne C$
then making inner blowing up at an edge of $C$ one can make its weight 0. If $C$ is a linear 0-vertex with neighbors of weight $w_1$ and $w_2$
then making an inner blowing up at an edge of $C$ and contracting $C$ we get a reconstruction $[[w_1,0,w_2]] \to [[w_1-1,0, w_2+1]]$.
Similarly, if $C$ is end $0$-vertex with a neighbor of weight $w$ one can get $[[0,w]] \to [[0, w+1]]$ or   $[[0,w]] \to [[0, w-1]]$. The last three
reconstructions around a $0$-vertex are called elementary transformations. The next straightforward fact will be useful.

\bprop\label{elementary}{\rm (see also \cite[Section 2]{FKZ})}. {\rm (1)} Let $C_1+C_2+C_3$ be a chain with weights $w_1, 0, w_2$. Then there exist elementary
transformations such that $C_1$ and $C_3$ 
are not blown down in this process and
one has the following change of weights $[[w_1,0,w_2]]\to [[w_1+w_2,0,0]]$.

{\rm (2)} Let $C_1+C_2 + \ldots +C_n$ be a chain (resp. a circular graph) with weights $0, 0, w_3, \ldots , w_n$. Then there exist  elementary transformations such that $C_i, \ldots ,
C_n$ are not blown down (where $3\leq i\leq n$)
 and one has the following change of weights $$[[0,0,w_3, \ldots ,w_n]]\to [[w_3, \ldots ,w_{i},0,0, w_{i+1}, \ldots , w_n]]$$
$$( {resp. \, \, \, } ((0,0,w_3, \ldots ,w_n))\to ((w_3, \ldots ,w_{i},0,0, w_{i+1}, \ldots , w_n)) \, ).$$ 
\eprop

  We say that a linear graph is a standard chain if it has one the following forms:
  $[[0_{2k+1}]]$, $[[w_1, \ldots , w_n]]$,
and $[[0_{2k}, w_1, \ldots , w_{n-1}]]$, where $n\geq 1$ $k\ge 0$ and every $w_i \leq -2$, see (cf. \cite{FKZ} Definition 2.13). Here the subscript for $0$ simply means the number of consecutive zeros.

We shall need later the following consequence of \cite[Theorems 2.15 and Theorem 3.1]{FKZ}.

\bprop\label{reconstruction} Let $\Gamma_i, \, i=1,2$ be minimal  non-circular graphs, $Br(\Gamma_i)$ be the set of
branch points of $\Gamma_i$, and $\Gamma_2$ admit a reconstruction from $\Gamma_1$. Then 

{\rm (i)} $Br(\Gamma_i)$ is an invariant of the reconstruction, i.e., none of the vertices of $Br(\Gamma_1)$ is blown down in this procedure
and they are transformed bijectively onto $Br(\Gamma_2)$ with preservation of valency;

{\rm (ii)}   there is a bijection between connected components of $\Gamma_1 \ominus Br(\Gamma_1)$
and of $\Gamma_2 \ominus Br(\Gamma_2)$ such that every connected component $\Gamma_2^0$ of $\Gamma_2 \ominus Br(\Gamma_2)$
is obtained from the corresponding connected component $\Gamma_1^0$ of $\Gamma_1 \ominus Br(\Gamma_1)$ by a sequence of blowing up and blowing down;

{\rm (iii)}  every minimal weighted graph  $\Gamma_1$ can be reconstructed into some minimal weighted graph  $\Gamma_2$ such that each connected component of $\Gamma_2 \ominus Br(\Gamma_2)$
is a standard  chain;

{\rm (iv)} if $\Gamma_1^0$ and $\Gamma_2^0$ are standard graphs then 
the reconstruction of $\Gamma_1^0$ into $\Gamma_2^0$ can be achieved by elementary transformations.

\eprop 

Since chains $[[w_1, \ldots , w_n]]$ with every $w_i \leq -2$ do not admit nontrivial elementary transformations we have the following.

\bcor\label{cor.reconstruction} Let $\Gamma_i, Br(\Gamma_i)$, and $\Gamma_i^0$ be as in Proposition \ref{reconstruction}.
Suppose that every weight in $\Gamma_1^0$ is at most -2.
Then any  relatively minimal reconstruction\footnote{That is, a reconstruction without unnecessary blowing up.}  of $\Gamma_1$ into $\Gamma_2$ induces
the identity transformation of $\Gamma_1^0$ into $\Gamma_2^0$. In particular, when every weight in
the graph $\Gamma_1 \ominus Br(\Gamma_1)$ is at most -2  any relatively minimal reconstruction is  the  identity transformation of $\Gamma_1$ into $\Gamma_2$.
\ecor


\section{$\PP^1$-fibrations} \label{P1fibrations}

The results of Section \ref{theorem1.3}  suggest that in order to classify complete algebraic vector fields on a semi-affine smooth surface $X$ we need
to understand
fibrations of $X$ with general fibers $\C$ or $\C^*$, for short $\C$- or $\C^*$-fibrations. They can be extended to $\PP^1$-fibrations (i.e., fibrations
with general fiber isomorphic to $\PP^1$) of a smooth completion $\bX$ of $X$. Hence, in this section we present some general results on $\PP^1$-fibrations.

\bnota\label{not2.1} For the rest of the section we suppose that $\bar f: \bX \to B$ is a $\PP^1$-fibration of
a smooth projective surface over a smooth complete curve.
Let $\bD$ be a connected curve in $\bX$  of simple normal crossing (SNC) type and $X=\bX \setminus \bD$  be semi-affine.
We suppose that $f=\bar f|_X$ is  a $\C$- or $\C^*$-fibration on $X.$  \enota

A classical result about $\PP^1$-fibrations is the following, see Proposition 4.3 in \cite{BPVdV}:

\bthm \label{ruled}
Let $\bX$ be a smooth compact surface and $C$ be a smooth rational curve in $\bX$. If $C^2 = 0$, then there exists a regular birational map $\varphi: X\rightarrow
Y$, where $Y$ is ruled (a $\PP^1$-bundle over a curve), such that $C$ meets no exceptional curve of $\varphi$, and $\varphi(C)$ is a general fiber of $Y$.
\ethm

In particular, this theorem states that singular fibers of a $\PP^1$-fibration are contractible to a rational curve and, thus, their dual graphs do not contain
cycles. The following lemma gives a slightly more precise statement about the singular fibers of a $\PP^1$-fibration.

\blem\label{irreducible} In the Setting \ref{not2.1} let $\Gamma$ be the dual graph of a fiber $F=\bar f^{-1}(b)$ for some $b \in B$.

{\rm (1)} 
Let  $E$ be a vertex of $\Gamma$ that is reduced
in $\bar f^* (b)$.  Then $\Gamma \ominus E$ is contractible  and, furthermore, after
this contraction the weight of $E$ becomes 0.

{\rm (2)} Let $E_1$ and $E_2$ be vertices of $\Gamma$ that are reduced in $\bar f^* (b)$ and let
$\Gamma^0$ be the smallest linear subgraph of $\Gamma$ containing $E_1$ and $E_2$.
Then $\Gamma \ominus \Gamma^0$ is contractible.

{\rm (3)} Let $C$ be a vertex of $\Gamma$ that has multiplicity 2
in $\bar f^* (b)$. Then $\Gamma$ can be contracted to a graph
of the form\\[2ex]
$$
\hspace{1.5cm}\bigskip \cou{C_n \, \, \, \, \, \, \,  }{ \, \, \, \, \, \, \, \, \, }\nlin\cshiftup{
E_{1}}{-2}\slin\cshiftdown{E_{2}}{-2}\lin\cou{C_{n-1}}{-2}\lin\ldots\lin\cou{C_{2}}{-2}\lin\cou{C_{1}}{-1} \hspace{0.5cm},\hspace{0.5cm} 
$$\\
where $C_1$ is the proper transform of $C$ and the weight of $C_n$ is $-2$ (resp. $-1$) for $n\geq 2$ (resp. $n=1$).

\elem

\bproof The first statement can be found, say, in \cite[Lemma 2.11.2]{Miy}. For the second statement we contract all components of
$\Gamma \ominus E_1$ that do not contain $E_2$ and vice versa. By abuse of notation denote the resulting graph by the same symbol $\Gamma$. Assume that $\Gamma \ominus \Gamma_0$ is still not empty
and it does not contain linear $(-1)$-vertices. It is enough to show that this is impossible. 

Since any connected component of $\Gamma \ominus E_1$ must be contractible by (1) this assumption implies that 
$E_1$ and similarly $E_2$ are end-vertices of $\Gamma$. Thus, if $\Gamma \ne \Gamma^0$
there exists a vertex $E\in \Gamma^0 \ominus (E_1 \cup E_2)$ such that it is a branch point of $\Gamma$. 
Since the graph is contractible to $E_1$
the component of $\Gamma \ominus E$ containing $E_2$ must be contractible to a point in the curve $E$. This implies that $E$ is reduced in  $\bar f^* (b)$
(indeed, if $E$ has multiplicity at least 2 so does $E_2$ because 
$E_2$ is obtained from a point in $E$ by a sequence of blow-ups). By (1) all components of $\Gamma\ominus E$ are contractible contradicting the assumption that
$\Gamma \ominus \Gamma^0$ does not contain linear $(-1)$ vertices which yields (2). 

 By Theorem \ref{ruled} contracting consequently $(-1)$-vertices one can transform $\Gamma$ into a reduced $0$-vertex.
In particular, during this contraction the proper transform of $C$ becomes a linear $(-1)$-vertex.  Thus, we can  suppose this is true
from the beginning. Consider two cases: (i) $C$ is not an end vertex of $\Gamma$ and (ii) $C$ is an end vertex of $\Gamma$.
Let $E_1$ and $E_2$ be the neighbors of $C$ in (i).  Note that $E_1$ and $E_2$ are reduced in $\bar f^*(b)$
since the multiplicity of $C$ in $\bar f^*(b)$ is the sum of the corresponding multiplicities of $E_1$ and $E_2$.
Hence, by (2) $\Gamma$ can be contracted to the chain $E_1+C + E_2$
and its weights must be $-2,-1,-2$ since otherwise $\Gamma$ cannot be contracted to a $0$-vertex. 
In (ii) contracting $C$ we observe that its neighbor must have multiplicity 2 in $\bar f^*(b)$. Thus, replacing $C$ by this neighbor and
using the induction by the number of vertices we conclude that $\Gamma$ can be contracted to a graph of the form presented in
statement (3) which concludes the proof.

\eproof

\bdefi If $f : X \to B$ is a $\C^*$-fibration then $\bD$ either contains two sections $B_1$ and $B_2$ of $\bar f$
(we call this case untwisted) or it contains a curve $B_0$ such that $\bar f |_{B_0}: B_0 \to B$ is a ramified double cover of rational curves
(so-called twisted case).
\edefi

\brem\label{lem.2singfib}
It is worth mentioning that in the twisted case  a fiber  $\bar f^{-1}(b), b \in B$ meets $B_0$ at one point if and only if $b$ is a ramification point
of the morphism $\bar f |_{B_0}: B_0 \to B$. 
The number of ramifications points is determined by
the Riemann-Hurwitz formula and when $B_0 \simeq \proj^1$ there are exactly two of them. 
\erem

\bprop\label{fiber} In the Setting \ref{not2.1} let $F=\bar f^{-1}(b)$ be a fiber contained in $\bD$.
Suppose that the dual graph $\Gamma$ of $\bD$ does not contain linear $(-1)$-vertices different from
irreducible curves on which the restriction of $\bar f$ is non-constant  (we call such $\bD$ pseudo-minimal). 
Let $\Gamma_0$ be the smallest subgraph of $\Gamma$ 
that contains all components of $F$ and their neighbors. 

{\rm (1)} If $f$ is a $\C$-fibration then $F$ is a $0$-vertex and $\Gamma_0$ is a linear chain $F+B'$ where $B'$
is a section of $\bar f$.

{\rm (2)} If $f$ is an untwisted $\C^*$-fibration then $F$ is a $0$-vertex and $\Gamma_0$ is a linear chain $B_1+F +B_2$.
 
{\rm (3)} Let $f$ be a twisted $\C^*$-fibration.

\hspace{.5cm} {\rm (3a)} Suppose that $F$ meets $B_0$ at two points. Then $F$ is a $0$-vertex
and $\Gamma_0$ is a cycle consisting of $F$, and $B_0$ joined by two edges. 

\hspace{.5cm} {\rm (3b)} Suppose that $F$ meets $B_0$ at one point. Then $F$ is a linear chain $C_1+E+C_2$ where $C_1$ and $C_2$
are $(-2)$-vertices, $E$ is a $(-1)$-vertex, and $\Gamma_0 \ominus E$ contains three components $C_1, C_2$, and $B_0$.
That is, the following
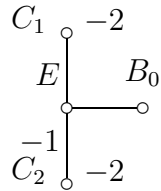
\vspace{1cm}

$$
\bigskip
\cou{E \, \, \, \, \, \, \,  }{-1 \, \, \, \, \, \, \,\, \,\, }\nlin\cshiftup{C_1}{-2}\slin\cshiftdown{C_2}{-2}\lin\cou{B_0}
\quad$$

\vspace{1cm}
\noindent is the form of $\Gamma_0$ (in the rest of the paper subgraphs of $\Gamma$ satisfying the assumptions of (3b) will be called subgraphs
of type  $\Gamma_*$).

\eprop

\bproof Note that in (1) $B'$ is the only neighbor of $F$ in $\Gamma$ since otherwise $f$ is not a $\C$-fibration. The component $E$ of $F$
that meets $B'$ is reduced in $\bar f^*(b)$ because $f^*(b) \cdot B'=1$. By Lemma \ref{irreducible} one can contract $F$ to $E$
but since $F$ does not have linear $(-1)$-vertices we have $F=E$.

Since $f$ is a $\C^*$-fibration in (2), by  similar reasons $B_1$ and $B_2$ are the only neighbors of $F$ in $\Gamma$ and the component $E_i$ of $F$
that meets $B_i$ is reduced in $\bar f^*(b)$. By Lemma \ref{irreducible} and pseudo-minimality $F$ must be a
linear chain with $E_1$ and $E_2$ being end vertices. In particular, $E_1$ and $E_2$ are linear vertices of $\Gamma$.
Note that if $E_1 \ne E_2$ then both of them are $(-1)$-curves since the curve $\overline{F \setminus E_i}$ is contractible by Lemma \ref{irreducible}.
This contradicts pseudo-minimality. Thus, $E_1=E_2=F$ which yields (2).

In (3a) as before we have $B_0$ as the only neighbor of $F$ in $\Gamma$ and the neighbors of $B_0$ in $F$ are $E_1$ and $E_2$
which are reduced in $\bar f^*(b)$ (indeed we have $f^*(b) \cdot B_0=2$). Again, by Lemma \ref{irreducible} and pseudo-minimality 
we see that $F=E_1=E_2$, i.e., being an irreducible fiber of $\proj^1$-fibration $F$ is the 0-vertex which is (3a).

In (3b) $B_0$ is again the only neighbor of $F$ in $\Gamma$. More precisely $B_0$ is a neighbor of a vertex $E$ in the graph $\Gamma_1$ of $F$
and $E$ has multiplicity 2 in $\bar f^*(b)$ since for the double section $B_0$ one has $\bar f^*(b) \cdot F=2$. Since $\Gamma_1$ is
contractible to a 0-vertex it contains a linear $(-1)$-vertex. The assumption on $(-1)$-vertices implies that this vertex is not in $\Gamma_1 \ominus E$.
Hence, $E$ is a linear $(-1)$-vertex in $\Gamma_1$. Let $C_1$ and $C_2$ be its neighbors in $\Gamma_1$. Note that they are reduced in $\bar f^* (b)$. 
Indeed, since $E$ is contractible to the points of intersection of the images
of $C_1$ and $C_2$ the multiplicity of $C$ is the sum of multiplicities of $C_1$ and $C_2$, i.e., 2=1+1. By Lemma \ref{irreducible} one can contract all
components of
$\Gamma_1 \ominus (C_1 \cup E \cup C_2)$ and we are done.
\eproof

\brem\label{3a} 
Let $C$ be a smooth rational curve in $\bX$ with $C^2=0$ (i.e., it is a fiber of a $\PP^1$-fibration by Theorem \ref{ruled}).  The following converse of Proposition \ref{fiber} is
true.

(1) Let $C$ be an end vertex of $\Gamma$, then $X$ admits a $\C$-fibration $f$ such that $C$ is a fiber of $\bar f$.

(2) Let $C$ be a linear vertex  of $\Gamma$ with two distinct neighbors $B_1$ and $B_2$. Then
$X$ admits an untwisted $\C^*$-fibration $f$ such that $C$ is a fiber of $\bar f$ and $B_1$ and $B_2$ are sections of $\bar f$.

(3a) Let $C$ be a linear vertex  of $\Gamma$ with one neighbor $B_0$ only (i.e., $C$ is joined with $B_0$ by two edges). Then $X$ admits a twisted
$\C^*$-fibration $f$ such that $C$ is a fiber of $\bar f$ and $B_0$ is the double section $\bar f$ which intersects $C$ transversally in two points.

(3b) Let  $C_1+C+C_2$ be a linear subgraph of $\Gamma$ as in (3b) with  $E$ replaced by $C$ and $B_0$ being the only neighbor of $C$ different from $C_1$ and $C_2$. Then $X$ admits a
twisted $\C^*$-fibration $f$ such that $C_1\cup C \cup C_2$ is a singular fiber of $\bar f$ and $B_0\subset \bD$ is the double section of $\bar f$.

\erem

We need one more technical fact  for Section \ref{homogen}.

\blem\label{Cfeather} In the Setting \ref{not2.1} let
$\bar f : \bX \to B$ be a pseudo-minimal extension of $f$, and $U$ be the union of components of $\bD$  on which the restriction of $\bar f$ is not constant.   
Suppose that $E_1$ and $E_2$ are the only irreducible components of
$\bar f^{-1}(b), \, b \in B$ that meet the divisor $U$ (where we allow equality $E_1=E_2$)
and that $\Gamma^0$ is as in Lemma \ref{irreducible}.

{\rm (i)} Let $F$ be
an irreducible affine\footnote{That is, a component that survives the Remmert reduction of $X$.} component of $f^{-1}(b)$ whose closure is not  a vertex in $\Gamma^0$.
Then $F\simeq \C$.

{\rm (ii)}  
If the dual graph $\Gamma'$ of $\bar f^{-1}(b)\cap \bD$ is not contractible
then $f^{-1} (b)$ is a singular fiber of the fibration $f$ (in other words, for every neighborhood $W$ of $b$ in $B$ the map $ f|_{ f^{-1} (W)} :
 f^{-1} (W) \to W$ is not a locally trivial fibration even in the sense of differential topology).  
\elem

\bproof    The dual graph $\Gamma''$ of the fiber $f^{-1} (b)$ is a tree since otherwise $\Gamma''$ cannot be contracted to a $0$-vertex.
 Let $\Gamma^1$ be the subgraph of $\Gamma''\setminus \Gamma^0$ consisting
of the vertices corresponding to the irreducible components of $\bar f^{-1}(b)$ meeting $f^{-1} (b)$. Since $X$ is semi-affine and, thus,
by the Lefschetz hyperplane section theorem $\bD$ is connected 
a vertex $\bF \in\Gamma^1$ can be only an endpoint of $\Gamma'' \ominus \Gamma^0$ which is not a neighbor of $\Gamma^0$. This implies 
that $\bF$ is an endpoint in $\Gamma''$ which yields (i).

If $f^{-1} (b)$ is not singular then this fiber is reduced irreducible. Its closure is can be viewed as a vertex $V$ in $\Gamma''$ such that
$\Gamma'=\Gamma''\ominus V$. By Lemma \ref{irreducible} $\Gamma'$ is contractible which is the desired conclusion.

\eproof


\section{On non-affine base of Riccati fibration}\label{non-affine}

The aim of this section is to show that under the assumptions of Theorem \ref{GR.theorem} with $B\simeq \PP^1$, there is another Riccati fibration for the same vector field which is not surjective.
The first step will be to prove that the $Z$ is a toric  surface.  The existence of a Riccati fibration with affine base in this case has been established in \cite{KLL} (heavily using Brunella's result for $\C^2$).

Let us prove  first several technical facts. 

\blem\label{new.l1} Let $X$ be a normal compact complex surface, $M$ be a finite subset of $X$, and $X^*=X\setminus M$.
Suppose that  $\nu$ is a complete holomorphic vector field on $X^*$. Then $\nu$ extends to a similar field on $X$.
\elem

\bproof Let $\varphi_t : X \to X, \, t \in \C$ be the  flow of $\nu$ and $x \in M$. When $t$ has a small absolute value then 
by continuity $\varphi_t$ transforms a punctured Stein neighborhood $U_x$ of $x$ into a similar Stein neighborhood $U_x'$.
Thus, by Hartogs' theorem $\varphi_t$ extends to $x$ which implies the desired conclusion.
\eproof

\blem\label{new.l2} Let the curve $B= g(Z)$ in  Theorem \ref{GR.theorem} be isomorphic to $\PP^1$ and let $\pi : X \to Z$ be a minimal resolution
of singularities of $Z$. Then there exists a smooth SNC-completion $\bX$ of $X$ so that the lift of $\sigma$ to $X$ extends to
an algebraic vector field $\bar \sigma$ on $\bX$.

\elem

\bproof By the assumption of Theorem \ref{GR.theorem} there exists a complete nontrivial algebraic vector field $\nu$ on $B$ for which
$f_* (\sigma) =\nu$. Let $O$ be the open orbit of $\nu$ in $B$ (i.e., either $O\simeq \C$ or $O \simeq \C^*$) and let $U =f^{-1}(O)$.
If $F=g^{-1} (z)$ for a general point $z \in O$ (i.e.,  either $F\simeq \C$ or $F \simeq \C^*$) 
then using the  flow of $\sigma$ one can see that $U$ is biholomorphic to $F\times O$ over $O$. In particular, all fibers of 
the morphism $g|_U$ are reduced
and isomorphic to $F$. Let us extend $g |_U$ to a proper morphism $\bar g: \bU \to O$ from a smooth surface $\bU$ which is, of course, a $\PP^1$-fibration. 
Since the closure of $F$ in $\bU$ is a reduced component of the fiber $\bar g^{-1}(z)$ we can contract all other components of this fiber
by Lemma \ref{irreducible} (1). That is, $\bU$ is naturally isomorphic to $\PP^1\times O$. 

Furthermore, one can see that because of completeness the 
field $\sigma$ is tangent to $\bU \setminus U$, i.e., it can be extended to a complete algebraic field on $\bU$. 
Thus, we have an algebraic extension of $\sigma$ to the surface $Z'=Z \cup \bU$. Since this extension vanishes at the singular points
we can blow these points up and obtain a complete algebraic vector field $\delta$ on the smooth surface $X'=X\cup \bU$ that extends $\sigma|_\bU$.
Extend the natural projection  $g' : X' \to B$ further to a proper morphism $\bar g' : \bX \to B$ such that it coincides with $\bar g$ over $O$.
Note that the curve $(\bar g')^{-1} (B \setminus O)\setminus X'$ has a negative definite intersection form by the Zariski lemma \cite[page
90]{BPVdV}. Thus, by the Grauert criterion \cite{BPVdV} it can be contracted. The result of this contraction $\rho : \bX \to \hX$ is a normal complex 
surface $\hX$ such that it contains $X'$
and $\hX \setminus X'$ is finite. By Lemma \ref{new.l1} $\delta$ extends to a complete holomorphic vector field on $\hX$.
Blowing up the singular points of $\hX$ we can lift this complete holomorphic vector field to a complete holomorphic vector field $\bar \sigma$
on $\bX$. Since $\bX$ is a complete algebraic surface this field $\bar \sigma$ is algebraic \cite{BPVdV} which concludes the proof.
\eproof

 \brem\label{non.aff.base.r1} The assumption that $B$ is a complete curve is crucial for Lemma \ref{new.l2}.
In its absence one can construct a complete vector field on a smooth affine surface $Z$ that do not admit a regular extension to any completion
$\bZ$ of $Z$. Indeed, for $Z=\C^* \times \C$ with a coordinate $x$ and $z$ on the first and the second factors respectively 
consider a vector field of the form $\nu= x \frac{\p}{\p x} - x^{-k}(P(x)z +Q(x))\frac{\p}{\p z}$ where $k\geq 0$ and $P$ and $Q$ are polynomials.
In this case $B=\C^*$ and $g: Z \to B$ is the natural projection. We can suppose that $g$ extends to $\bar g : \bZ \to \PP^1$.
However, it is not difficult to see that the extension of $\nu$ has poles either on $\bar g^{-1} (0)$ or on $\bar g^{-1} (\infty)$.
\erem

\blem\label{new.l3} {\rm (1)} Let $\bar g : \bX \to B$, $X$ and $\bar \sigma$ be as in Lemma \ref{new.l2} (in particular, $\bar \sigma$ does not have a rational
first integral and $\bar g$ is a surjective map onto $B\simeq \PP^1$). 
Suppose that $Y=\PP^1 \times \PP^1$ where the first (resp. second) factor is equipped with an affine coordinate $z$ (resp. $w$).
Then there exists a non-regular birational map $\chi : \bX -\to Y$ such that it transforms $\bar \sigma$ into a complete holomorphic vector field $\mu = \mu_1 + \mu_2$ on $Y$
where either 

(a) $\mu_1=\alpha z \frac{\p}{\p z}$ and $\mu_2=\beta w \frac{\p}{\p w }$ with $\alpha$ and $\beta\in \C$ being  linearly independent over $\Q$, or

(b) $\mu_1= \frac{\p}{\p z}$ and $\mu_2=\gamma w\frac{\p}{\p w}$ with $\gamma \ne 0$.

{\rm (2)} Furthermore, let $C_1=\{ w =0 \}$, $C_2=\{ z =\infty \}$, $C_3=\{ w =\infty \}$, and $C_4=\{ z=0  \}$ be the lines in $Y$. 
Let $R=\bigcup_{i=1}^4C_i$ in case (a) and  $R=\bigcup_{i=1}^3C_i$ in case (b). Then the restriction of $\chi$ over $Y\setminus R$
is an isomorphism and $\chi^{-1} : Y \dashrightarrow \bX$ is regular outside the set $P$ of singular points of $R$.

{\rm (3)} Suppose that $D=\bX \setminus X$ and $T=D\cup \bar g^{-1} (B\setminus O)$ (where $O$ is as in the proof of Lemma \ref{new.l2}). 
 Then $\chi (T) \subset R$. 

\elem

\bproof By \cite[Section 6, Theorem 6 (iii)]{Bru00}  there exists $\chi : \bX -\to Y$ such that it transforms $\bar \sigma$ into a complete holomorphic vector field $\mu = \mu_1 + \mu_2$ on $Y=\PP^1\times \PP^1$ where $\mu_i$ is the lift of a complete vector field on the $i$-th factor of $Y$. This implies that, say, $\mu_1$ is of the form
$l(z)  \frac{\p}{\p z}$ where $l(z)$ is a quotient of  linear polynomials. However, some of these fields $\mu$ have rational first integrals (e.g., $\mu = \alpha \frac{\p}{\p z}+ \beta \frac{\p}{\p w}$ with
$\alpha$ and $\beta \in \C$ has the rational first integral $\beta z-\alpha w$). Disregarding such fields 
 we see that up to automorphisms of the factors $\PP^1$ we are left with the choices presented in cases (a) and (b) which is (1).
 
 It follows from the proof of \cite[Section 6, Theorem 6 (iii)]{Bru00} that $\chi$ is a composition of a sequence of monoidal transformations
 and contractions of $(-1)$-curves where each monoidal transformation occurs at a point where the vector field (induced by) $\bar \sigma$ is zero
 while for every contracted $(-1)$-curve such a field either is tangent to it or vanishes on it (in particular, the image of this curve after the contraction
 is a point where the consequent induced vector field is zero).
 
 Note that $\mu$ does not vanish on $Y\setminus R$ and none of integral curves of $\mu$ contained in $Y \setminus R$ is algebraic, i.e.,
such curves  cannot be contracted under $\chi^{-1} : Y \to \bX$. Thus, $\chi$ is an isomorphism over $Y\setminus R$. Furthermore, $\mu$ vanishes only
at $P$ and, therefore, the indeterminacy points of $\chi^{-1}$ are contained in $P$ which concludes the proof of (2).

By construction $T$ consists of irreducible components invariant under the foliation induced by $\bar \sigma$. Hence, $\chi (T)$ consists of irreducible components
invariant under the foliation induced by  $\mu$.
Furthermore, $\chi (T)$ is an algebraic curve since $D$ is. Since all algebraic curves tangent to $\mu$ are contained in $R$ we get (3).
\eproof

\bdefi\label{new.d1}  Let $\nu$ be a  holomorphic vector field on a smooth complex surface $X$ vanishing at a point $x_0 \in X$.
Then $x_0$ is a reduced singularity of $\nu$ if the linear part of $\nu$
at $x_0$ has eigenvalues $\lambda_1, \lambda_2$ such that either they are nonzero and $\lambda_1/\lambda_2 \notin  \Q_+$ or
$\lambda_1 \ne 0=\lambda_2$ (e.g., every point $p \in P$ from Lemma \ref{new.l3} is a reduced singularity for the vector fields $\mu$ described in that lemma).

\edefi

The main property of reduced singularities we are going to exploit is the following. Let $\pi : \tX \to X$ be the blowing up of $X$ at $x_0$
and $E \simeq \PP^1$ be its exceptional divisor. Then $\nu$ induces the vector field $\tilde \nu$ on $\tX$ such that it is tangent to $E$,
vanishes exactly at two points of $E$, and has both of these singularities reduced.\footnote{One can extract this from the following computation.
Suppose that $\alpha z \frac{\partial}{\partial z} + \beta w \frac{\partial}{\partial w}$ is the linear part of $\nu$ in a local coordinate system $(z,w)$ at $x_0$.
Then we can consider a local system $(\xi, \eta)$ on $\tX$ such that $(z,w)=(\xi, \xi \eta)$. The restriction of $\tilde \nu$ to $E$
is given by $\alpha \xi \frac{\partial}{\partial \xi} + (\beta -\alpha ) 
\eta \frac{\partial}{\partial \eta}$. }

\blem\label{new.l4} Let the assumption of Lemma \ref{new.l3} hold. Then there exist a smooth surface $W$ and morphisms
$\varphi : W \to \bX$ and $\psi : W \to Y$ such that $\chi = \psi \circ \varphi^{-1}$ and the dual graph of the curve $\psi^{-1}(R)$ is either circular or linear.
Furthermore,  the vector fields $z \frac{\p}{\p z}$ and $w \frac{\p}{\p w }$ on $Y$ induce holomorphic vector field $\delta_1$ and $\delta_2$ on $W$ such
that they are tangent to every component of $\psi^{-1}(R)$.
\elem

\bproof The existence of $\varphi : W \to \bX$ and $\psi : W \to Y$ such that $\chi = \psi \circ \varphi^{-1}$ follows from the Zariski theorem \cite{BPVdV}.
We suppose also that these morphisms are chosen so that there are no curves in $W$ that are contracted by both $\varphi$ and $\psi$.
Then according to the proof of \cite[Section 6, Theorem 6 (iii)]{Bru00} $\psi$ is a composition of a sequence of monoidal transformations 
each of which occurs at a point where the vector field (induced by) $\mu$ vanishes. Since $\mu$ vanishes at $P$ only
the first of these transformations $Y'\to Y$ takes place at a point $p \in P$. It transforms $R$ to a curve $R'$. Since the dual graph $\Gamma$ of $R$ is linear or circular 
and $p$ is a double point of $R$ the dual graph $\Gamma'$ of $R'$ is obtained from $\Gamma$ by an inner transformation (see Section \ref{dualgraphs}
for the definition of an inner transformation)
and, therefore, $\Gamma'$ remains linear or circular. Let $\mu'$ (resp. $\delta_1'$, resp. $\delta_2'$) be the vector field on $Y'$ induced by $\mu$
(resp. $z \frac{\p}{\p z}$,  resp. $w \frac{\p}{\p w }$) and $E \subset R'$ be the exceptional curve of
the transformation. As we mentioned before there are only two points on $E$ where $\mu'$ vanishes (and they are reduced singularities of $\mu'$).
The same is true for $\delta_1'$ and $\delta_2'$.
Since the whole curve $R'$ is invariant under the  flow of $\mu'$ (resp. $\delta_1'$, resp. $\delta_2'$) these two points must be automatically the points where $E$ meets
other irreducible components of $R'$. That is, they are again double points of the simple normal crossing curve $R'$. As we mentioned before
the consequent monoidal transformation must occur at a zero of $\mu'$  (and such a zero is simultaneously a zero of $\delta_1'$ and of $\delta_2'$).
This leads to an inner transformation of $\Gamma'$.
In particular, the absence of branch points for the dual graphs is preserved on this step and similarly on all consequent steps
while the lifts of $z \frac{\p}{\p z}$ and $w \frac{\p}{\p w }$ remain holomorphic and tangent to the preimage of $R$.
This yields the desired conclusion.
\eproof

\bprop\label{new.p1} Let the assumptions of Theorem \ref{GR.theorem} hold and the fibration $g$ being onto $B\simeq \PP^1$. 
Then $Z$ is a toric surface. 
\eprop

\bproof
By Lemma \ref{new.l3} we have $D'=\psi^{-1}(D)$  and $T'=\psi^{-1}(T)$ contained in $R'=\psi^{-1}(R)$.
By Lemma \ref{new.l4} the dual graph of $R'$ is either circular or linear. If $D' = R'$ then $T'\setminus D' = \emptyset$, i.e., $\bar g$ has no
fibers invariant under the foliation induced by $\bar \sigma$. Hence, $\bar g$ is not surjective contrary to the assumption. Thus,
the dual graph of $D'$ is linear.  Therefore, the dual graph of $D$ is linear which implies  already that $Z$ is a Gizatullin surface \cite{Gi}.
Furthermore, by Lemma \ref{new.l3} the restriction of $\chi$ over $Y \setminus R$ is an isomorphism which implies that
the curves contracted by the birational morphism $\varphi : W \to \bX$ are contained in $R'$. Since these curves are tangent to the vector fields
$\delta_1$ and $\delta_2$ on $W$ (from Lemma \ref{new.l4}) we can push these fields down to holomorphic vector fields on $\bX$ tangent to $D$.
That is, these resulting fields are induced (via $\chi$) by the vector fields $z \frac{\p}{\p z}$ and $w \frac{\p}{\p w }$ on $Y$ whose  flows commute
and generate the action of the torus $(\C^*)^2$ on $Y$. This implies that we have the torus action on $\bX$ and $X$. Every complete curve in $X$
has a negative self-intersection and, thus, it is invariant under the action. Contraction of such curves leads to a non-degenerate torus action on $Z$
which concludes the proof. 
\eproof

\bthm\label{KKtheorem}
 Let $\sigma$ be a  complete algebraic vector field on a normal affine surface $Z$ such that $\sigma$ does not
admit a rational first integral. Then there is a regular function $g: Z \to \C$  whose general fibers are isomorphic to $\C$ or $\C^*$ so that
the  flow of $\sigma$ sends fibers of $g$ to fibers of $g$. 
\ethm

\bproof
We have to consider the case when the fibration provided by Theorem \ref{GR.theorem} is surjective onto $B\simeq \PP^1$. By Proposition \ref{new.p1} this can only happen when we are dealing with a toric surface
$Z$. For toric surfaces the existence of another Riccati fibration with  affine base has been deduced in \cite{KLL} (Lemma 4.3)
from Brunella's original work. \eproof

\section{Rational first integral}\label{rational}

 The aim of this section is to classify generalized Gizatullin surfaces admitting nonzero complete algebraic vector fields that have  non-regular first integrals.

\bprop\label{field} Let $B$ be a germ of a smooth curve at point $o$ 
and  let  $\mu$ be the vector field on $B \times \PP^1_x$ given by $\mu=x{\partial}/{\partial x}$ where $\PP^1_x = \C \cup \infty$ and 
$\C$ is equipped with a coordinate $x$ (i.e.,  
the set of its zeros of $\mu$ 
is the union of $B\times\lbrace\infty\rbrace$
and $B\times\lbrace 0 \rbrace$). Suppose that $\pr : \tX \to B$ is a smooth $\PP^1$-fibration and $\psi : \tX \to B \times \PP^1_x$ is
a birational morphism over $B$ whose restriction over $B \setminus o$ is an isomorphism. Let the dual graph of $\pr^{-1}(o)$
be linear with end vertices meeting the proper transforms $\tB_1$ and $\tB_2$ of $B\times\lbrace\infty\rbrace$
and $B\times\lbrace 0 \rbrace$ respectively. Then $\mu$ induces a regular complete vector field $\tilde \mu$ on $\tX$ tangent to the fibers of $\pr$ and such that
its restriction to $\pr^{-1}(o)$ vanishes only at double points of the curve $\tB_1 \cup \pr^{-1}(o) \cup \tB_2$.
\eprop

\bproof We use induction on the number $k$ of irreducible components of $\pr^{-1}(o)$. If $k=1$ then $\psi$ is an isomorphism
and there is nothing to prove. Note that contracting a $(-1)$-curve in $\pr^{-1}(o)$ we obtain a birational morphism
$\sigma : \tX \to \breve{X}$ over $B$ such that $\psi$ factors through it and 
the dual graph of $\breve{\pr}^{-1}(o)$ (for the natural projection $\breve{\pr} : \breve{X} \to B$) is
linear with end vertices meeting the proper transforms $\breve{B}_1$ and $\breve{B}_2$ of $B\times\lbrace\infty\rbrace$
and $B\times\lbrace 0 \rbrace$ respectively. By assumption $\mu$ induces a regular complete vector field $\breve{\mu}$ on $\breve{X}$
tangent to the fiber of $\breve{\pr}$ and such that
its restriction to $\breve{\pr}^{-1}(o)$ vanishes only at double points of the curve $\breve{B}_1 \cup \breve{\pr}^{-1}(o) \cup \breve{B}_2$.

Note that $\sigma$ is a monoidal transformation with center also at one of these points.  Hence, $\breve{\mu}$  generates a field $\tilde {\mu}$ 
tangent to the fibers of $ \pr$ which is complete by construction. Its  flow preserves the curve $\tB_1 \cup \pr^{-1}(o) \cup \tB_2$
and, in particular, it  keeps the set $S$ of the double points of this curve fixed, i.e $\tilde \mu$ vanishes at these points. For any component $F$ of $\pr^{-1}(o)$
on which $\tilde \mu$ is not identically zero, $\tilde \mu$ does not vanish on $F \setminus S \simeq \C^*$ since no rational curve but $\C$ or $\C^*$ can
be an integral curve of a complete vector field. Therefore, it remains to show that $\tilde \mu$ has only isolated zeros on $\pr^{-1}(o)$.

By induction one can suppose that in some local coordinate system $(z,w)$ near a double point of $\breve{B}_1 \cup \breve{\pr}^{-1}(o) \cup \breve{B}_2$
the local equation of $\breve{B}_1 \cup \breve{\pr}^{-1}(o) \cup \breve{B}_2$ is $zw=0$ and
the field $\breve{\mu}$
coincides with  $$nz\frac{\partial}{\partial z} -m w\frac{\partial}{\partial w}$$  for natural $n$ and $m$ with $n+m\geq 1$. 
Furthermore, for a local coordinate
system $(\xi, \eta)$ on $\tX$ the map $\sigma$ is given by $(z,w)=(\xi, \xi \eta)$. The direct computation shows that
the local form of $\tilde \mu$ is
$$n\xi \frac{\partial}{\partial \xi} -(n+m) \eta \frac{\partial}{\partial \eta}$$
which implies the desired conclusion.
\eproof

\brem\label{semi} Note that $\mu$ is semi-simple, i.e., its  flow is an algebraic $\C^*$-action on $B \times \PP^1$.
Hence, the  flow of $\tilde \mu$ is also an  algebraic $\C^*$-action on $\tX$ and $\tilde \mu$ is semi-simple as well.
\erem

\bnota\label{nota.rational}
Let $\nu$ be a complete vector field on a smooth semi-affine surface $X$ with a rational first integral. Blowing up $X$  at the points
of indeterminacy (note that $\nu$ vanishes at such points) we can suppose that this integral is a regular morphism  $f:X\to B$ where $B$ is a complete curve.
In particular, $f$ is either a $\C$- or a $\C^*$-fibration.  Let $\hat f : \hX \to B$ be an extension of $f$ to a $\proj^1$-fibration on a smooth completion $\hX$ of $X$
by an SNC-divisor $\hD$ which is assumed to be pseudo-minimal. Similarly, we suppose that the union $\hE$ of complete curves in $\hX$ does not contain $(-1)$-curves tangent to $\nu$
(and we call such $\hE$ pseudo-minimal).
The extension of $\nu$ to $\hX$ will be denoted by $\hat \nu$ (i.e., $\hat \nu$ may have poles). We note that unless $f$ is a twisted $\C^*$-fibration 
the set of zeros of $\hat \nu$ contains \\

\noindent (a) either only one section $B_0$ of $\hat f$ (i.e a general integral curve of $\hat \nu$ is isomorphic to $\C$) or\\

\noindent (b)  two sections $B_1$ and $B_2$ of $\hat f$ (i.e a general integral curve of $\hat \nu$ is isomorphic to $\C^*$) \\

\noindent where the second option is automatic for
the untwisted $\C^*$-fibration.  Furthermore,
$B_0$ (resp. at least $B_1$) must be contained in $\hD$ since otherwise $X$ is not semi-affine.

\enota

\blem\label{lem.rational}
In the Setting \ref{nota.rational}  let either (a) or (b) hold. 
Suppose that  $\mu={\partial}/{\partial x}$ (resp. $\mu=x{\partial}/{\partial x}$) is the vector field on $B\times \PP^1_x$ in case (a) (resp. (b)).
Then there exists a birational map $\varphi : \hX \dasharrow B \times \PP^1$ over $B$ (in particular the restriction of $\varphi$ over some open Zariski
dense subset $B^*$ of $B$ is an isomorphism) for which

{\rm (1)} $\mu$ induces a rational vector field $\hat \mu$ on $\hX$ such that for some rational function $p$ on $B$ one has $\hat \nu =\hat f^*(p) \hat \mu$;

{\rm (2)}  in case (b) the restriction of $\hat \mu$ to any fiber of $\hat f$ has a finite
number of zeros;

{\rm (3)}  if in case (b) $\hat \mu|_{X} $ is regular  then $\hat \mu$
is semi-simple and the dual graph of the curve $B_1 \cup \hat f^{-1}(b) \cup B_2$ is linear for every $b \in f(X)$.

\elem
\bproof 
By Theorem \ref{ruled} there is a birational morphism $\tau : \hX \to Y$ onto a ruled surface $Y$ over $B$.
Recall that $Y$ is a locally trivial $\PP^1$-fibration over $B$ \cite[Prop. V.2.2]{Har}, i.e., there is a Zariski dense open
subset $B^*$ of $B$ such that the preimage $Y^*$ of $B^*$ under the projection $Y \to B$ is naturally isomorphic
to $B^* \times \PP^1$ (more precisely, such $B^*$ can be chosen as a neighborhood of any $b_0$ for which the fiber $f^{-1}(b_0)$ is
not a singular one). Furthermore, reducing $B^*$ one can suppose that the proper transforms of $B_1$ and $B_2$  in $Y^*$
(in case (b)) are disjoint and the restriction of $\tau$ to $\hX^* =\tau^{-1}(Y^*)$ yields an isomorphism $X^* \simeq Y^*$. 
Using the freedom of choice of this isomorphism $Y^* \simeq B^* \times \PP^1$
one can suppose now that these proper transforms are $B^*\times\lbrace\infty\rbrace$
and $B^*\times\lbrace 0 \rbrace$ respectively (and in case (a) the proper transform of $B_0$ is $B^*\times\lbrace\infty\rbrace$).
Extending  the isomorphism  $Y^* \simeq B^* \times \PP^1$ we get $\varphi : \hX \dashrightarrow B \times \PP^1$. Let us show that
$\varphi$ is the desired rational map.

By the Zariski theorem there is a surface $W$ that dominates both $\hX$ and $B \times \PP^1$ over $B$. In particular,
in case (b) the proper transforms $B_1'$ and $B_2'$ of $B_1$ and $B_2$ in  $W$ are disjoint (since $B\times\lbrace\infty\rbrace$
and $B\times\lbrace 0 \rbrace$ are).
By Lemma \ref{irreducible}(2)
the dual graph of every singular fiber of the natural morphism $\kappa : W \to B$ is contractible to the minimal linear subgraph joining two vertices
meeting $B_1'$ and $B_2'$ respectively. Making all such contractions one gets a morphism $W \to \tX$. It remains to note that
the morphism $W \to B \times \PP^1$ must factor through $W \to \tX$ since one wants to keep the proper transforms of $B_1'$ and $B_2'$ disjoint.
That is, $\varphi =\psi \circ \chi$ where $\psi : \tX \to B\times \PP^1$ is a birational morphism over $B$ and $\chi : \hX \dasharrow \tX$  is the birational map (that factors through $W$). 

The image of $\hat \nu$ under isomorphism $\hX^* \simeq Y^*$ yields a complete vector field $\nu_0$ on $Y^*\setminus B^*\times\lbrace\infty\rbrace \simeq B^* \times \C_x$
tangent to the fibers of the natural projection onto $B^*$.
Hence, in case (a) the restriction of the field $\nu_0$ to every fiber must be proportional
to $\mu={\partial}/{\partial x}$  and, thus, there is a regular function $p$ on $B^*$ for which $\nu_0=p\mu$. The extension $p$ to $B$ is the desired
rational function. 
In case (b) the argument about existence of $p$ is similar but with $\mu= x{\partial}/{\partial x}$ which yields (1).

By Proposition \ref{field} and Remark \ref{semi} in case (b) $\mu$ induces a semi-simple field $\tilde \mu$ on $\tX$ whose restriction on any fiber $\pr^{-1}(b)$
of the natural projection $\pr : \tX \to B$ has zeros only at the set $Z_b$ of double points of the curve $\tB_1 \cup \pr^{-1}(b) \cup \tB_2$ where $\tB_i$ is the proper transform of $B_i$.
By construction the image $I$ of the exceptional divisor $F$ of the morphism $W \to \tX$ does not meet $Z_b$ for any $b\in B$.  
That is, $\tilde \mu$ does not vanish at any point of $I$
and, therefore, $\tilde \mu$ induces a rational vector field on $W$ with poles on $F$.
In particular,  the restriction of this rational vector field  to any fiber of  $\kappa : W\to B$ has a finite number of zeros which yields (2).

Assume that for some $b\in f(X)$ the fiber $\kappa^{-1}(b)$ has a nonlinear dual graph, 
or, equivalently, $F \cap \kappa^{-1}(b)$ is not empty. The proper transform of $F \cap \kappa^{-1}(b)$ in $\hX$ is not contained in $\hD$
because of pseudo-minimality assumption.  That is, there is a component $C$ in $F \cap \kappa^{-1}(b)$
whose proper transform $\hat C$ in $\hX$ meets $X$. Hence, $\hat \mu |_{X}$ is not regular because of poles on $\hat C$.
Thus, for regularity one needs $\chi$ to be an isomorphism which yields (3).
\eproof

We need to consider two essentially different cases: when $f: X \to B$ is surjective and when it is not. 
As the next claim shows the assumption that there is a vector field $\nu$ tangent to fibers of a $\C$- and $\C^*$-fibration in Setting \ref{nota.rational}
is automatically fulfilled in the non-surjective case.

\bprop\label{tangent.field} Let $f: X \to B$ be a non-surjective morphism from a smooth semi-affine surface $X$ into $B$ with  general fibers isomorphic to $\C$ or to $\C^*$.
Then there is a complete algebraic vector field tangent to the fibers of $f$. Furthermore, one can choose this field so that it does not vanish on a given general fiber
$E=f^{-1}(b_0)$ of $f$.
\eprop

\bproof Suppose first that $f$ is an untwisted $\C^*$-fibration (resp. $\C$-fibration). As we showed in the proof of Lemma \ref{lem.rational} for every regular value $b_0 \in B$ of 
$f$ there is a Zariski neighborhood $B^* \subset f(X) \subset  B$ for which $Z=f^{-1}(B^*)$ is naturally isomorphic to $B^* \times \C^*$ (resp. $B^* \times \C$) over $B^*$. 
Vector field $\hat \mu$ from Lemma \ref{lem.rational} may have poles on $X \setminus Z$. However, since
$f(X)$ is affine one can choose a regular
function $h$ on $f(X)$ with prescribed orders of zeros at points of $f(X) \setminus B^*$ so that $h\hat \mu$ yields a regular complete
algebraic vector field on $X$ tangent to the fibers of $f$. Choosing $h$ so that $h(b_0) \ne 0$ we get the second statement.

In the twisted case consider a proper extension $\hat f : \hX \to B $ of $f$ to an SNC-completion $\hX = X \cup \hD$ of $X$. Then $\hD$ has the only one irreducible
component $B_0$ on which the restriction of $\hat f$ is not a constant. Furthermore, $p:=\hat f|_{B_0}$ makes  $B_0$ a ramified double cover of $ B$.
Hence, one has a $\Z_2$-action $\alpha$ on $B_0$ for which $B \simeq B_0/\alpha$. Let $X_0 =X \times_B B_0$ and $q: X_0 \to B_0$, $r: X _0 \to X$ be the induced
morphisms. Since $r$ makes $X_0$ a ramified double cover of $X$ we have again a $\Z_2$-action $\beta$ on $X_0$ for which $X=X_0/\beta$.
Consider an $\alpha$-invariant 
Zariski dense open subset $B_0^* \subset p^{-1}(f(X))\subset  B_0$
for which $q^{-1}(B_0^*)$ is naturally isomorphic to $B_0^* \times \C^*$ where the second factor is equipped with a coordinate $x$.
Then the restriction of $\beta$ to $q^{-1}(B_0^*)$ 
is given by $\beta (b,x)=(\alpha (b), g(b,x))$ where $(b,x) \in B_0^*\times \C^*$ and for a fixed $b$ the function $g(b,x)$ is a coordinate on  $\C^*$.
That is, $g(b,x)$ coincides either with $e(b)x$ or with $e(b)/x$ with $e(b)$ being a non vanishing regular function on $B_0^*$.
However, the first possibility must be disregarded because we deal with a twisted case.

In order to construct a complete algebraic vector field on $X$ tangent to the fibers of $f$ on $X$ it suffices to construct a complete
algebraic vector field on $(p\circ q)^{-1}(f(X))$ tangent to the fibers of $q$ and such that its restriction to $q^{-1}(B_0^*)$
is $\beta$-invariant. Note that the field $x\frac{\partial}{\partial x}$ is mapped to  $-x\frac{\partial}{\partial x}$ under automorphism $x \to e(b)/x$ of $\C^*$.
Hence, for every function regular function $h$ on $B_0^*$ that is $\alpha$-antisymmetric  the field $h x\frac{\partial}{\partial x}$ is invariant
with respect to the $\beta$-action. Choosing this function $h$ on $B_0$ so that its extension to the affine curve $p^{-1}(f(X))$ has zeros of sufficiently high order at points
of $p^{-1}(f(X))\setminus B_0^*$ we guarantee the regular extension of this field to  $(p\circ q)^{-1}(f(X))$ which yields the first statement.
For the second statement choose $B_0^*$ so that it contains $p^{-1}(b_0)$ and require that $h$ does not vanish on this set. Hence, we are done.  
\eproof

In particular, the presence of $\C$- and $\C^*$-fibrations in  Theorem \ref{KKtheorem}  implies the following.

\bcor\label{quasi.flexible}
Every normal semi-affine surface $X$ with a complete algebraic vector field which does not possess  a rational first integral has automatically
an open orbit.
\ecor

\blem\label{lem.sur.rational}
Let $\hat \nu$ and $\hat f: \hat X \rightarrow B$ be as in Setting \ref{nota.rational} and $f : X \to B$ be surjective.
Then 


{\rm (i)} $f$ is an untwisted $\C^*$-fibration;

{\rm (ii)}  every singular fiber of $\hat f$ has a linear dual graph $C_1 + \ldots + C_n + C + C'_{n'} + \ldots + C'_1$ where $n,n' \geq
1$, $C$ is the only
$(-1)$-vertex of this fiber (i.e., $C_i\cdot C_i$ and   $C_i'\cdot C_i' \leq -2 $), $C_i \subset \hat D$ and $C'_i \subset \hat E$ while $C$ meets both $X$ and $\hD$;

{\rm (iii)}   
the field $\hat \nu$ (and, therefore, $\nu$) is semi-simple.
\elem

\bproof
Assume that contrary to (i) we deal with case (a) from Setting \ref{nota.rational}. Let $\varphi: \hat X\dashrightarrow B \times \PP^1$, $\hat \mu$ and $p$ be as in Lemma
\ref{lem.rational}, and let $B^*$ and $\hX^*\simeq B^* \times \PP^1$ be as in the proof of that Lemma. In particular, $p$ is regular on
$B^*$ and $\hat \nu =\hat f^*(p)\hat \mu$
is locally nilpotent on $\hX^* \setminus B_0 \simeq B^*\times \C_x$ (since $\hat \nu =p \frac{\partial} {\partial x}$).
That is, the  flow of this complete vector field $\hat \nu$ induces an algebraic $\C_+$-action on $\hX^* \setminus B_0$ and, therefore, on $X=\hX \setminus \hD$.
However, the algebraic quotient $X//\C^+$ must be affine (since $X$ is a semi-affine surface) while the surjectivity of $f|_{X}$
implies that this quotient is $B$. Hence, since $B$ is complete
we have to disregard case (a).

Suppose now that we are in case (b) and use  the notation from Lemma \ref{lem.rational}.  By surjectivity of
$f\vert_{X}$ every singular fiber $\hat f^{-1} (b)$ contains an irreducible component $C$ that meets $\hD$ but is not contained in $\hD$. 
By Lemma \ref{lem.rational} $\hat \mu$ does not vanish identically on $C$.  Hence, $p$ has no pole in $b$
since the field  $\hat\nu =\hat f^*( p) \hat \mu^*\vert_{X}$ is regular. The absence of poles on the complete curve $B$
means that $p$ is constant and, thus, we can suppose that $\hat \mu=\hat \nu$.

This implies in particular that $\hat \mu$ is regular on $\hat X \setminus \hat D$. By Lemma \ref{lem.rational}(3)  
the fibers of $\hat f$ are linear chains of rational curves. Note also that $\hD \cup \hE$ is the union of $B_1$, $B_2$, and all irreducible components
of singular fibers  with exception of components similar to $C$. Therefore, $B_1$ and $B_2$ belong to different connected components of $\hD \cup \hE$.
Recall that $B_1\subset \hD$.
Since $\hD$ is connected by the Lefschetz theorem, we see that the connected component containing $B_1$ (resp. $B_2$) coincides with $\hD$ (resp. $\hE$).

Assume that $n=0$, i.e., $C$ meets section $B_1$ which implies that $C$ is irreducible in the fiber $\hat f^*(b)$. By Lemma \ref{irreducible} and pseudo-minimality
of $\hE$ one has $n'=0$ contrary to the fact that $\hat f^{-1}(b)$ is singular. Thus, $n\geq 1$ and similarly $n'\geq 1$ which yields (ii).

It remains to exclude the case of twisted $\C^*$-fibration, i.e the case when $\hat D$ contains a double section $B_0$. Let us replace $\hat f : \hX \to B$ by
the natural morphism  $\hX \times_B B_0\to B_0$ and also replace $\hat \nu$, $X$, and 
$\hD$ by their lifts to $\hat X \times_B B_0$. Then two sections of the modified morphism $\hat f : \hX \to B$ are contained  in $\hD$ contrary to the argument before,
i.e., we have (i) and (ii).
 
Lemma \ref{lem.rational} (3) and the equality $\hat \nu = \hat \mu$ imply (iii) which concludes the proof.
\eproof

\brem\label{FZ} (1)
It follows from the proof that $\hE$ is connected. Thus, the Remmert reduction $X_0$ of
$X$ has only one singularity which is automatically a fixed point of an elliptic $\C^*$-action associated with $\nu$. Hence, Lemma \ref{lem.sur.rational} 
can be also obtained from the description of normal $\C^*$-singularities according to \cite{OW,Pink}. It can be extracted as well from
the DPD-presentation for $\C^*$-surfaces due to Flenner and Zaidenberg \cite{FZ}.  

(2) Furthermore, if $X_0$ is smooth then the Luna
slice theorem implies that $X_0 \simeq \C^2$ \cite{Lu}.

(3) In the case when $X_0$ is  not smooth 
consider the field $\nu_0$ induced by $\nu$ on $X_0$ and the rational first integral $f_0$ induced by $f$. 
One can see that the surjectivity of $f : X \to B$ is equivalent to the fact that $f_0$ is not regular on $X_0$.
\erem

\bnota\label{second.fibration} Suppose we are in  Setting \ref{nota.rational}   and let 
$g: X \to B'$ be another $\C$- or $\C^*$-fibration on $X$ over a complete curve $B'$
such that a general fiber of $g$ is not contained in
a general fiber of $f$ and vice versa.
Let $\hat g : \hX' \to B'$ be a proper extension of $g$ to an SNC-completion
of $X$ by a pseudo-minimal divisor $\hD'$.
\enota

\bprop\label{prop.2sur.rational} Assume the Setting \ref{second.fibration}. Then $B$ (and, therefore, $B'$) is a rational
curve and $X$ is a rational surface. Furthermore,  if $f: X \to B$ is surjective then either

{\rm (1)}  the dual graph of $\hD$ (and, therefore, of $\hE$) is linear, or

{\rm (2)} the dual graph of $\hD$ is of form

\vspace{1cm}

$ \hspace{1.5cm}\bigskip \cou{B_1 \, \, \, \, \, \, \,  }{-1 \, \, \, \, \, \, \, \, \, }\nlin\cshiftup{\tilde
C_1}{-2}\slin\cshiftdown{\tilde C_2}{-2}\lin\cou{C_0}{w_0}\lin\cou{C_1}{w_1}\lin\ldots\lin\cou{C_{n}}{w_{n}} \hspace{0.5cm},\hspace{0.5cm} 
 $
 \vspace{1cm} 
 
\noindent where $w_i\leq -2$ for  $i\geq 0$.  In particular, in (2) $\hat f$ has three singular fibers $\hat f^{-1}(b_0), \hat f^{-1}(b_1)$, and $ \hat f^{-1}(b_2)$
containing $C_0, \tilde C_1$, and $\tilde C_2$ respectively.

Moreover, $g$ is not surjective.
\eprop
\bproof

By assumption the restriction of $f|_F : F \to B$ to a general fiber $F$ of $g$ is not constant. Hence, $B$ is rational since $F$ is rational.
This implies also that $X$ is rational because $f$ is a  $\C$- or $\C^*$-fibration.

Assume that $f$ is surjective and   $\hat f$ has three or more singular fibers or equivalently that the dual graph of $\hat D$ is not linear.
By Lemma \ref{lem.sur.rational} $B_1$ is the only branch point  of this graph and the weight of any other vertex is at  most
$-2$. The identical automorphism of $X$ extends to a rational map
$\hat X \dashrightarrow \hat X'$, i.e., there is a reconstruction of the dual graph $\Gamma$ of $\hD$ into a dual graph $\Gamma'$ of
$\hD'$.

By pseudo-minimality the weights of linear vertices of $\Gamma'$
are at most $-2$. Hence, Corollary \ref{cor.reconstruction} implies that  if $\Gamma'$ is minimal then it coincides with $\Gamma$ and otherwise
it is obtained from $\Gamma$ by a sequence of inner and outer blowing up.
In particular, $\Gamma'$ cannot contain a linear $0$-vertex and if it has a branch vertex of weight $-1$ then it is a proper transform of $B_1$.

However, if $\hat g : X \to B'$ is not surjective then by Proposition \ref{fiber} there must be either  a subgraph of type $\Gamma_*$ from Proposition \ref{fiber} (3b)    or a $0$-vertex in $\Gamma'$.
In combination with the previous argument Proposition \ref{fiber}  leads to the graph in (2).

Thus, it remains to consider the case when $g : X' \to B'$ is surjective.
Let us derive a contradiction by showing that the fibrations $\hat f$ and $\hat g$ must coincide in this case.
As before  we see that $\Gamma'$ is obtained from $\Gamma$ by a sequence of blowing-ups. Hence, by  Corollary \ref{cor.reconstruction}   of $\hD$ we can suppose
that $\hX = \hX'$.

Let $F$ (resp. $G$) be a general fiber of $\hat f$ (resp. $\hat g$). 
It suffices to show that
$F$ is equivalent to $G$ in $H^2( \hat X , \Z )$ (because the linear systems $|F|$ and $|G|$ induce $\hat f$ and $\hat g$ respectively).
Let  $S$ (resp. $S'$) be  the elements of this cohomology group corresponding to the vertices of  the dual graphs of $\hat E \setminus B_2$  (resp. $\hat D
\setminus B_1$).
Since $\hat X$ is contractible to a surface ruled
over $B$ we see that the elements of $S$ and $S'$ together with $F$ and $B_2$ form
a basis of $H^2( \hat X, \Z )$. In particular,  in this cohomology group $G =kF+lB_2 + M+M'$ where $M$ (resp. $M'$) is an integer linear combination
of elements from $S$ (resp. $S'$). Note that the restriction of the intersection form to $S$ (resp. $S'$) is negative definite by the Zariski lemma \cite[page
90]{BPVdV}.
By the same reason $F \cdot C_i = F \cdot C_j' = 
G \cdot C_i = G \cdot C_j' =0$. Hence, $(M')^2 = G \cdot M'=0$ which implies that $M'$ is the zero divisor. Since $G \cdot B_1=1$ we have
$k=1$ and $G =F + lB_2+M$.  Then equalities $B_2 \cdot G  =1, G \cdot M=0$ and $(G )^2=0$ imply that 
$$l (B_2)^2+ B_2 \cdot M =0, \, \, M^2 + lB_2 \cdot M =0, \, \, \, {\rm and} \, \, \, M^2 +l^2 (B_2)^2+ 2lB_2 \cdot M +2l =0.$$ 
In turn the last three equalities yield $l=0$ and, therefore, $M^2=0$. The semi-negativity implies that $M=0$ and $G =F$ which is
the desired conclusion.
\eproof

\bprop\label{toric} Let us be in Setting \ref{second.fibration} and $f: X \to B$ be surjective. Suppose also that $\hD$ has a dual graph different from the graph
in Proposition \ref{prop.2sur.rational} (2). Then the Remmert reduction $X_0$ of $X$ 
is a normal toric surface.
\eprop

\bproof  Since $\hD$ is the boundary divisor of a semi-affine surface the Nakai-Moishezon criterion implies that $B_1^2\geq -1$ and, furthermore, if $B_1^2=-1$ then
$\hD$ can be contracted to a curve which contains an irreducible component of zero weight and some other components of negative weight.
 If $B_1^2 >0$ then after several inner blowing-ups we can make its weight equal to zero.
 In both cases (and also in the case of $B_1^2=0$) one of the end vertices of the dual graph of the resulting curve $\tD$ is of negative weight and it corresponds to the end-vertex $C_n$ of $\hD$ which is
 automatically contained
 in some singular fiber $\hat f^{-1}(b)$. Let $C_1 + \ldots + C_n + C + C'_{n'} + \ldots + C'_1$ be the dual graph of $\hat f^{-1}(b)$ as in Lemma \ref{lem.sur.rational}. 
 Then the dual graph of $\tD$ is of form $F_1+F_2+ \ldots + F_k+ C_n $ and furthermore, using elementary transformations
 from Proposition \ref{elementary} we can suppose that $F_1$ is of zero weight. In particular, $F_1$ generates a $\PP^1$-fibration $\tilde f : \tX \to \PP^1$ 
 on the resulting surface $\tX$. Note that by construction $F_2$ (which may be equal to $C_n$) is a section of this fibration
 and the only singular fiber  $\tilde f^{-1} (a)$  of $\tilde f$ is
 the union $G$ of $\hE$ and the components $F_3,\ldots ,F_k, C_n, $  and $C$. Indeed $G$ is connected and does not meet $F_1$ (i.e., it is contained in
 a fiber of $\tilde f$), and
 there is no complete irreducible curve different from a component of $G$ that
 does not meet both $F_1$ and  $F_2$ (which would be the case if $G$ is not the whole fiber or in the case of another singular fiber). Furthermore, since $X =\hX \setminus \hD =\tX \setminus \tD$ and
 the dual graph of $\hD$ (or $\tD$) is linear the Remmert reduction $X_0$ of $X$ is a normal Gizatullin surface and in terminology of \cite{FKZ.a} the curve $\tD \cup \tilde f^{-1}(a)$
 is called the extended boundary divisor. A normal Gizatullin surface (say $X_0$) is toric if and only if the dual graph of the extended boundary divisor for its minimal resolutions of singularities  
 (say $X$) is linear
\cite[Lemma 2.20]{FKZ.a}. This is exactly the case of the dual graph of $\tD \cup \tilde f^{-1}(a)$ and we are done.
\eproof

\bprop\label{nontoric} Let us be in Setting \ref{second.fibration}  and $f: X \to B$ be surjective. Suppose also that the dual graph of $\hD$ is as
in Proposition \ref{prop.2sur.rational} (2). Then 

{\rm (1)} $X$ is a surface with an open orbit;

{\rm (2)} the  flow of $\nu$ preserves the fibers of $g$;

{\rm (3)} the Remmert reduction $X_0$ of $X$ is not generalized Gizatullin unless $n=0$ and $w_0=-2$ in which case the dual graph of the
curve $\hE $ has the following form

\vspace{1cm}

$ \hspace{4.5cm}\cou{C_0'}{-2\,\,\,\,\, }\lin\cou{ \, \, \, \, \, \, \,\, \,\,B_2}{\,\,\,\,\,\,\, \,\,\, -2}
\nlin\cshiftup{\tilde C_1'}{-2}\slin\cshiftdown{\tilde C_2'}{-2}
$
 
 \vspace{1cm} 
\noindent whence the dual graph of every singular fiber $ \hat f^{-1}(b_i)$ coincides with $[[-2,-1,-2]]$ for $i=0,1,2$.
Furthermore, such a surface $X$ is unique up to isomorphism.

\eprop

\bproof By Proposition \ref{tangent.field} there is a complete algebraic vector field $\mu$ tangent to the fibers of $g$.
Hence, $\nu$ is not proportional to $\mu$ which yields (1).

By Lemma \ref{lem.sur.rational} $\hat \nu$ (resp. $\nu$) is a semi-simple field on $\hX$ (resp. $X$). In particular, the  flow $\Phi_t$ of $\hat \nu$
preserves the curve $F_\infty =\tilde C_1 \cup B_1 \cup \tilde C_2$ which is the fiber of $\hat g$. Thus, for any time parameter $t$ and every 
other fiber $F$ of $\hat g$ the image $\Phi_t(F)$ does not meet $F_\infty$ which implies that $\hat g$ is constant on $\Phi_t(F)$, i.e., it is
again a fiber of $\hat g$. This yields (2).

By Lemma \ref{lem.sur.rational}  the dual graph of $\hat f^{-1}(b_1)$ is $$\tilde C_1 + \tilde C + \tilde C_{n'}''+ \ldots + \tilde C_1''$$ where $\tilde C$ is the
only $(-1)$-vertex, any other weight is at most $-2$, and $\tilde C_1^2=-2$. Since this fiber of a $\PP^1$-fibration must be
contractible to a 0-vertex, contracting $\tilde C$ and $\tilde C_1$ consequently we see that $n'=1$ and $(C_1'')^2=-2$. The same is
valid for  the fiber $\hat f^{-1}(b_2)$, and in the case of $n=0$ and $w_0=-2$ it is also true for $\hat f^{-1}(b_0)$.

Assume that $w_0 \leq -3$. Consider the dual graph $$C_0+ C_1+ \ldots +C_n+C +C_{n'}'+ \ldots +C_1'$$ of $\hat f^{-1}(b_0)$.  Then the connected curve 
$$G=C_1+ \ldots +C_n+C +C_{n'}'+ \ldots +C_1' + B_2+ \tilde C_1' + \tilde C_2'$$ is a fiber of $\hat g$
since it contains all complete curves in $\hX$ that do not meet the fiber $F_\infty$ or the section $C_0$. Note that $B_2$ is branch point of the dual graph of $G$
and all branches but $$\cB =C_1+ \ldots +C_n+C +C_{n'}'+ \ldots +C_1'$$
are non-contractible (because we know already that the weights of $ \tilde C_1' $ and $\tilde C_2'$ are $-2$). 
Thus, the latter must be contractible since any fiber is contractible to a $0$-vertex and furthermore, after this contraction
the weight of $B_2$ must become $-1$ in the graph $[[-2,-1,-2]]$.
Hence, as in the proof of Proposition \ref{fiber} we conclude that $B_2$ is a multiple component of $\hat g$ and, therefore, each component of $\cB$   
is also multiple. In particular, $C$ is a multiple component of $g$, i.e., $C\cap X$ is a singular fiber of $g$.

Suppose that $\eta$ is a complete algebraic vector field without a rational first integral. By Theorem \ref{KKtheorem} 
there must a $\C^*$- or $\C$-fibration over $\C$ such that its fibers are transformed into each other by the  flow of $\eta$ and $g$ is the only candidate for this fibration
(since its minimal extension $\hat g$ to a $\PP^1$-fibration is determined uniquely by the form of $\hD$).
However, by  Lemma \ref{lem6.5} 
the curve $C\cap X$ is preserved by this flow.

Similarly, if $C\cap X$ is a component of a fiber of a rational first integral of some complete vector field it is preserved by the  flow of this field
(in particular this is true for $\nu$ and $\mu$).  If we have $w_0\leq -3$ or $n\geq 1$ then $g$ is the unique $\C^*$- or $\C$-fibration on $X$
and, thus, the curve $C \cap X$ is preserved by all complete vector fields. Since the image of $C\cap X$ in $X_0$ is a curve $X_0$ cannot be generalized Gizatullin when $w_0\leq -3$. 
or $n\geq 1$. 

Thus, we have the desired form of the dual graph of $\hat E$ except for the fact that the weight $k$ of $B_2$ is still unknown. Note first that
$k\leq -2$ since otherwise the intersection matrix of $\hE$ is not negative definite contrary to the Grauert
criterion of contractibility \cite[Theorem III.2.1]{BPVdV}.  Then contraction of
fibers $\hat f^{-1}(b_0), \hat f^{-1}(b_1)$, and $ \hat f^{-1}(b_2)$ to $0$-vertices transforms $\hX$ into a Hirzebruch surface $\hX'$ with proper transforms $B_1'$ and $B_2'$ of $B_1$ and $B_2$
as disjoint sections. Their weights $k_1$ and $k_2$ depends on the choice of contraction of these three fibers but in any case $k_1+k_2=-1+k+3=k+2\leq 0$.
Since these sections are disjoint it follows from
\cite[page 518]{GrHa} that one of them is of negative weight (say $-m$) and the weight of the other is at least $m$ which yields $k=-2$ and we are done
with the form of the graph of $E$.

Furthermore, we see now that $\hX'$ can be chosen as $\PP^1\times \PP^1$ and $\hX$ is obtained from $\hX'$ by some standard blowing up in three fibers of a natural morphism $ \hX' \to \PP^1$.
Since the group of automorphisms of $\PP^1$ acts transitively on the triples of distinct points we get the claim about uniqueness. 
\eproof

\bexa\label{D4}  Let $X_0$ be from Proposition \ref{nontoric}. Then its singularity (whose minimal
resolution is given by the graph of $\hE$) is a du Val singularity of  type $D_4$
where a singularity of type $D_{n+1}$ is locally isomorphic to the hypersurface $yx^2+y^n+z^2=0$ in $\C^3_{x,y,z}$. In fact  for $n=3$ one can see that $X_0$ is globally isomorphic to
the hypersurface $y(x^2+y^2)+z^2=0$. The elliptic $\C^*$-action on it is given by 
$(x,y,z) \to (\lambda^2x,\lambda ^2y,\lambda^3z)$ for $\lambda \in \C^*$ and it corresponds to a semi-simple field $\sigma$.
Three $\C^*$-fibrations associated with the three different strings $[[-2,-1,-2]]$ in the graph of $\hD$ are given by
the functions $y, x+\sqrt{-1}y$, and $ x-\sqrt{-1}y$. By Proposition \ref{tangent.field} there are complete algebraic vector fields $\sigma_1,\sigma_2, \sigma_3$ on the hypersurface tangent to
the fibers of these functions respectively. Note that the only curve tangent to both $\sigma$ and $\sigma_1$ is given by $y=0$, i.e this
curve is invariant under their  flows. However, it is not invariant under the  flows of
$\sigma_2$ or $\sigma_3$ and, hence, one can check that the natural action of $\AAutH (X_0)$ has the smooth part of $X_0$ as an open orbit.
In particular, $X_0$ is generalized Gizatullin (another proof of this fact will be considered in the last section) and, hence, it is the surface described in
Proposition \ref{nontoric} (3). For $n\geq 4$ there is only
one regular $\C^*$-fibration on the corresponding surface (given by function $y$) which is, therefore, not a generalized Gizatullin surface, only a surface with an open orbit.
\eexa

\bthm\label{thm.tor} Let $X_0$ be a  generalized Gizatullin surface such that for a complete algebraic vector field $\nu_0$ on $X_0$
there is a surjective rational first integral $f_0: X_0 \dasharrow B$ into a complete curve $B$. Then 

{\rm (1)} either $X_0$ is toric (and in 
particular a Gizatullin surface) or $X_0$ is the Remmert reduction of $X$ from Proposition \ref{nontoric};

{\rm (2)} up to a constant nonzero factor $\nu_0$ is semi-simple.

\ethm

\bproof Statement (2) follows from \ref{lem.sur.rational}(iii). For (1) consider a birational morphism $X \to X_0$ from a semi-affine surface $X$ such that $f_0$ induces a surjective morphism $f : X \to B$ with
general fibers isomorphic to $\C$ or $\C^*$.
If a complete algebraic vector field $\mu_0$ on $X_0$ non-proportional to $\nu_0$ has also  a rational first integral 
we can suppose that it induces a similar morphism $g : X \to B'$. Then Propositions \ref{toric} and \ref{nontoric} imply the desired conclusion.

If $\mu_0$ does not have a rational first integral then by  Theorem \ref{KKtheorem}  there is a morphism $g_0: X_0\to B'$
with general fibers isomorphic to $\C$ or $\C^*$. Since $f_0$ has indeterminacy points its general fibers are different from the general fibers of $g_0$
and, therefore, the fibration $g: X \to B'$ (induced by $g_0$) is different from $f: X \to B$. By Proposition \ref{tangent.field} there is a complete vector
field tangent to the fibers of $g$ and we are done again by Proposition \ref{toric}.
\eproof

Now we are in the position to prove our second important result from the introduction, which in the case of $\C^2$ can be extracted from \cite{Bru04}.

\bproof (of Theorem \ref{mthm1i}) Assume that $X$ has no open orbit, i.e., $\nu$ must have a rational first integral $f: X \dashrightarrow B$ by Proposition \ref{tangent.field}.
Then we have (1).

For (2)  consider a surface $X$ with an open orbit. Because of Proposition \ref{tangent.field} we can suppose again that $\nu$ has 
a rational first integral.  If such an integral is surjective then by Theorem \ref{thm.tor} we may deal
either with a toric surface or with the surface from Proposition  \ref{nontoric}. In the latter case the second statement of Proposition \ref{nontoric} yields the desired regular function $f: X \rightarrow \C$.
In the former case the existence of such $f$ was established in \cite{KLL} or it can be extracted from explicit description of algebraic $\C^*$-actions on Gizatullin surfaces 
in \cite{FKZ.a}. Now suppose that
every complete algebraic vector field on $X$ has a regular rational first integral. Then such integrals for non-proportional fields lead
to different $\PP^1$-fibrations on a completion of $X$ satisfying the assumption of Proposition \ref{prop.2sur.rational} whence $X$ is rational.
Hence, for any $\C$- or $\C^*$-fibration $X \to B$ one has $B \simeq \C$.
In any case this regular rational first integral for $\nu$ can be viewed as a desired $f: X \rightarrow \C$ in (2) and we are done.
\eproof


 \section{Dual graphs of curves in rational surfaces}\label{existinggraphs}

 The aim of this section is the following.

\bprop\label{dug.p0}  {\rm (1)} Let $\Gamma$ be a graph as in Theorem \ref{mthm} (3) (resp. (5)) such that $w_0\geq 0$.
Then $\Gamma$ is the dual graph of an SNC curve (with rational irreducible components) contained in a smooth complete rational surface 
if and only if the assumption ($\alpha$) (resp. ($\beta$)) from Theorem \ref{} holds.

{\rm (2)}   Let $\Gamma$ be a graph as in Theorem \ref{mthm} (2a).
Then $\Gamma$ is the dual graph of an SNC curve (with rational irreducible components) contained in a smooth complete rational surface 
if and only if the assumption ($\gamma$) from Theorem \ref{} holds.
 \eprop
  
 The proof requires some preparations.
  
  \blem\label{dug.l1}    
  Let $C_0$ be a smooth rational curve in a smooth complete rational surface $\bX$
  admitting a $\PP^1$-fibration $\bar f : \bX \to \PP^1$  such that the map $\bar f|_{C_0} : C_0\to \PP^1$ has
 degree 2. Then there is a regular birational map $\iota: \bX \to \hX$
  to Hirzebruch surface $\hX$
  such that its restriction  $\iota|_{C_0}:  C_0\to\hC_0:=\iota (C_0)$ is an isomorphism and $\bar f$ factors
  through a natural projection $\hat f : \hX \to \PP^1$. 
  \elem
   
\bproof  By the Riemann-Hurwitz formula there are  only two points $b_0, b_1 \in \PP^1$ such that $F_i=f^{-1} (b_i), \, i=0,1$
meets $C_0$ at one point.  Since the degree of $\bar f|_C : C \to \PP^1$ is 2
every other fiber $F$ of $\bar f$ meets $C_0$  transversely at two points and
the components of $F$ that meet $C_0$ are reduced in $\bar f^*(b)$
(there at most two of them, say $E_1$ and $E_2$). By Lemma \ref{irreducible} we can contract $F$  first to
a linear chain whose end points are $E_1$ and $E_2$ and then 
contract it further to $E_1$. Note
that this procedure does not produce the singularities on the image of $C_0$ since we contract consequently
$(-1)$-curves which are transversal to (the proper transform of) $C_0$ and which meet it at most at one point. By Theorem \ref{ruled}
we can contract also $F_0$ and $F_1$ to $0$-vertices which yields the desired Hirzebruch surface.
By construction $\iota|_{C_0}:  C_0\to\hC_0$ is injective and it remains to show that $\hC_0$ is smooth over $b_0$ and $b_1$.

If more than one irreducible component of $F_0$ meets $C_0$ then we can argue as in the case of a general fiber $F$.
Thus, we can suppose that only one irreducible component of $F_0$ meets $C_0$. If it is not transversal to $C_0$ then its multiplicity in $\bar f^* (b_0)$ is 1
and we can contract the rest of the fiber $F_0$ to it without affecting $C_0$. Thus, we can suppose that this component
meets $C_0$ transversely and its multiplicity in $\bar f^*(b_0)$ is 2. Then again by Lemma \ref{irreducible} (3) we can
contract $F_0$ to a chain $\tC_1+E+\tC_2=[[-2,-1,-2]]$ without producing a singularity and with $E$ meeting the proper transform of $C_0$ transversely.
Note that contracting consequently $E$ and $\tC_1$  we do not produce a singularity either by the same transversality argument.
One can deal with  $F_1$ in the similar manner
which concludes the proof.
\eproof
    
\bprop\label{dug.p1} Let $C$ be a smooth rational curve in a Hirzebruch surface $\hX=\F_n$
such that $\kappa|_C : C \to \PP^1$ has degree 2 where $\kappa : \hX \to \PP^1$ is the natural projection.
Then $C\cdot C=4$.
\eprop 

\bproof Treating closed curves in $\hX$ as elements of the N\'eron-Severi group  of $\hX$  we recall that this group
is generated by a fiber $\hF$ of $\kappa$ and the negative section $S$ of $\kappa$, i.e., $\hF^2=0, \hF\cdot S=1$
and $S^2=-n$. Since $\hF$ is a smooth rational curve the adjunction formula yields $(K_\hX +\hF)\cdot \hF=-2$ (where $K_\hX$
is the canonical divisor). Hence, $K_\hX =-2S+m\hF$ for some $m \in \Z$. On the other hand,
by the same reason $(K_\hX + S)\cdot S =-2$ which implies that $m=-2-n$ and 
$K_\hX =-2S-(n+2)\hF$. Since $\kappa|_C$ is of degree 2 we have $C=2S+l\hF$ for some $l \in \Z$. Since $C$ is also smooth rational
we have $(K_\hX + C)\cdot C =-2$ which implies that $ 2(l-n-2)=-2$. Hence, $l=n+1$.
Now we have $C^2=-4n+ 4(n+1)=4$ which concludes the proof.
\eproof
 
\bcor\label{dug.c1} Let $\Gamma_0$ be a graph as in Proposition \ref{fiber} (3b) and $\Gamma$ be a graph containing
$\Gamma_0$ such that $B_0$ is the only neighbor of the chain $C_1+E+C_2$ in $\Gamma$.
Suppose that $\Gamma$ is the dual graph of a SNC-curve in a smooth complete rational surface $\bX$ such
that all irreducible components of this curve are rational. Then the weight of $B_0$ is at most 2 (and it is exactly 2 when $\Gamma=\Gamma_0$).
\ecor  

\bproof The chain generates a $\PP^1$-fibration $\bar f : \bX \to \PP^1$. 
By Lemma \ref{dug.l1} there is a contraction $\iota : \bX \to \hX$ to a Hirzebruch surface $\hX$
in which the image $C$ of $B_0$ has self intersection $w=4$ by Proposition \ref{dug.p1}.
Since $\iota$ involves contraction of at least two consequent neighbors of $C$ (say, $E$ and $C_1$) we see that the weight of $B_0$
is at most $4-2=2$ which yields the desired conclusion.
\eproof

\brem\label{dug.r1}  For $\hat X=\F_n$ with $n \geq 2$ a curve $C$ as in Proposition \ref{dug.p1}  does not exist. 
Indeed, as we saw in the proof $C =2S + (n+1)F$ and, hence, $C \cdot S = 1-n<0$ which is absurd.
   For the quadric $\PP^1\times \PP^1$ such curve $C$ obviously exists (it is enough to consider the standard parabola
 in $\C^2\subset \PP^1\times\PP^1$ and it completion in the quadric).
  Choose a general line $S'=b_1\times \PP^1$
and a line $F=\PP^1\times b_2$ such that  $F\cap C\cap S$ is a point $a$.
 Blowing up the quadric at $a$ and contracting the proper transform of $F$ we get a Hirzebruch surface $\F_1$ for
 which the proper transform of $S'$ is the negative section. Note that the proper transform
 of $C$ is still smooth and its natural map to  $\PP^1$ is still of degree 2, i.e., we have also the desired curve for $\F_1$.
 \erem
 
 \bprop\label{dug.p3} Statement (1) of Proposition \ref{dug.p0} is true.
 
 \eprop
 \bproof Assume  $\bX$ is a smooth complete rational surface  containing a connected SNC-curve $G$ with rational components
 such that its dual graph $\Gamma$ is as  in Theorem \ref{mthm}  (3) or (5).  Then $\Gamma$ contains a subgraph
 $\Gamma_0$ as in Proposition \ref{fiber} (3) which contains in turn the chain $\tC_1+E+\tC_2$ generating 
a $\PP^1$-fibration $\bar f : \bX \to \PP^1$. The vertex $C_0$ of $\Gamma_0$ yields a smooth rational curve
such that $\bar f|_{C_0} : C_0 \to \PP^1$ has degree 2.  Making contractions in all fibers of $\bar f$ but the one 
containing the chain we obtain another $\PP^1$-fibration $\kappa : \hX \to \PP^1$  with all fibers but one (say, $F_0$)
irreducible and by Lemma \ref{dug.l1} we can suppose that the image $B_0$ of $C_0$ is smooth. That is,
the curve $F_0\cup B_0$ has the same dual graph as in Proposition \ref{fiber} (3) with the weight of $B_0$ equal to 2 (by Corollary \ref{dug.c1}).
 Note that  Remark  \ref{dug.r1}  implies that such a fibration
$\kappa$ and the curve $F_0\cup B_0$ can be constructed starting from a quadric.  

The preimage of $F_0 \cup B_0$ in $\bX$ must be a curve with a dual graph $\Gamma'$ such that
 $\Gamma$ is a subgraph of $\Gamma'$ with contractible  $\Gamma'\ominus \Gamma_0$ and
 $C_0\in \Gamma_0 \subset \Gamma$ being a linear vertex of $\Gamma$ with a nonnegative weight.
 Clearly if we can find such a graph  $\Gamma'$, a desired surface $\bX$ together with a curve $G$ exist.
 
 If $\Gamma'=\Gamma_0$ then we have configuration (3) from Theorem \ref{mthm}
 with $n=0$ and $w_0=2$. If $\Gamma=\Gamma_0$ and $\Gamma'\ominus \Gamma$ consists of
 one (resp. two) $(-1)$-vertices then we have the same configuration with $n=0$ and $w_0=1$ (resp. $w_0=0$). 
 In the case of $n\geq 1$ in configuration (3) $\Gamma\ominus \Gamma_0$ should be a linear chain $C_1+ \ldots + C_n$ with weights $w_i \leq -2$.
 To show that $\Gamma'$ containing such $\Gamma$ exists let $\Gamma'\ominus \Gamma$ be a disjoint collection of $(-1)$-vertices
 such that $C_i$ has $-w_i-2$ neighbors from this collection for $1\leq i \leq n-1$, while $C_n$ has $-w_n-1$ such neighbors.
 Then contraction of these $(-1)$-vertices transforms the chain into $[[-2,-2, \ldots, -2, -1]]$. Hence  $\Gamma'\ominus \Gamma_0$
 is contractible and this contraction yields $w_0=1$. Joining, if necessary, $C_0$ with an extra $(-1)$-vertex from $\Gamma'\ominus \Gamma$
 one can make $w_0=0$ which yields all possibilities for configuration (3) listed in  Theorem \ref{mthm} ($\alpha$).
 For configuration (5) $\Gamma\ominus \Gamma_0$ is of the form\\[2ex]
 $$
 \hspace{1.5cm}\bigskip
\cou{C_1}{w_1}\lin\ldots\lin\cou{C_{n}}{w_n}\lin\cou{ \, \, \, \, \, \, \, E'}{ \, \, \, \, \, \, \, k'}\nlin\cshiftup{\tilde
C_1'}{-2}\slin\cshiftdown{\tilde C_2'}{-2}
$$\\
where $n\geq 0$ and $w_i\leq -2$ for $i \geq 1$. Consider $n \geq 1$. To get a contractible $\Gamma'\ominus \Gamma_0$
one must have a contractible branch of $\Gamma'$ at $E'$ containing, say, $\tC_1'$ (which takes place if the branch consists of $\tC_1'$ and a $(-1)$-vertex).
Contraction of this branch yields a linear chain in which the proper transforms of $E'$, $\tC_2'$ and $C_n$ will be
denoted by $\breve E'$, $\breve C_2'$ and $\breve C_n$ respectively. Note that the weight of $\breve E'$ is $k'+1$. 
If $k'\leq -3$ then we deal with a linear chain as in configuration (3) and the same argument produces a contractible $\Gamma'\ominus \Gamma_0$
and $w_0=1$ or $w_0=0$.  If $k'=-2$ then $\breve E'$ is a $(-1)$-vertex and we can contract $\breve E'$ and $\breve C_2'$.
Note that this contraction changes the weight of $\breve C_n$ to $2+w_n$ and for further contraction we have to require that
the latter number is negative, i.e. $w_n\leq -3$. This yields all configurations in Theorem \ref{mthm} (5) with $n\geq 1$.

For $n=0$ consider the second fiber $F_1$ of $\kappa$ that meets $B_0$ at one point only.
Make two monoidal transformations of $\hX$ at $F_1\cap B_0$ and
the infinitely near point. Then the graph of the preimage of $F_0\cup B_0\cup F_1$ yields configuration
(5) with $n=0$ and $k'=-1$. Using additional blowing up we can also make $k'\leq-2$ which exhausts all configurations 
(5) listed in  Theorem \ref{mthm} ($\beta$) and, thus, concludes the proof.
 \eproof

\bprop\label{dg.p2} Statement (2) of Proposition \ref{dug.p0} is true.
\eprop

\bproof Consider the cycle $((0,0,0,0))$ which is the dual graph of an SNC-curve in the quadric whose components are rational.
This implies the part ``if" of the statement. 

Let $\Gamma=C_0 + C_1 +\ldots + C_n+C_{n+1}$ be a cyclic graph with $n \geq 2$, weights $w_0=w_{n+1}=0$ and $w_i \leq -2$ for $1\leq i\leq n$
such that  $\Gamma$ is the dual graph of an SNC curve
in a smooth complete surface $\bX$ with rational irreducible components.
By Theorem \ref{ruled} $C_0$ generates a $\PP^1$-fibration with $C_1$ and $C_{n+1}$ being sections of this fibration.
Furthermore, $C_2+ \ldots +C_n$ is contained in a fiber $F$ of the fibration. Note that the multiplicity of $C_n$ in $F$ is 1 since
$C_n$ meets the section $C_{n+1}$. By Lemma \ref{irreducible} $F$ is contractible to a zero vertex that is the proper
transform of $C_n$. This contraction leads to a curve with a circular dual graph $C_0'+C_1'+C_n'+C_{n+1}'$ where $C_i'$ is
the proper transform of $C_i$ and the weights of $C_0', C_n'$, and $C_{n+1}'$ are zero. Note that $C_1'$ is contained in a fiber of the fibration induced by $C_{n+1}'$.
Repeating the previous argument we can suppose that it coincides with this fiber, i.e its weight is zero which yields the desired conclusion.
\eproof

Now  Propositions \ref{dug.p3} and \ref{dg.p2} yield Proposition \ref{dug.p0}.

\bexa\label{dg.e1} (1) Let $n\geq 5$ and $k$ be the integer part of $n/2$. Let $\Gamma_0 =((0,0,w_1, \ldots , w_n))$  be a circular weighted graph with
$w_k$ and $w_{k+1}\leq -2$ and $w_i\leq -3$ for $i \ne k, k+1$. Then $\Gamma_0$ is a dual graph of an SNC-curve with rational components contained in a smooth rational surface.

Indeed, consider $n-2$ consequent inner blowing up in the cycle $((0,0,0,0))$ so that for an even (resp. odd) $i$
the $i$-th blowing up occurs at the left (resp. right) edge adjoint to the vertex created in the previous step. Then we get
the following cycle $((0,0, -3, \ldots , -3, -2, -1, -3 , \ldots , -3))$ where $-2$ and $-1$ occupy positions $k$ and $k+1$ or vice versa.
Now after appropriate outer blowing up this cycle can be transformed into the desired one.

(2) It is easy to show that the cycle $((0,0,-2,-2,-2,-2,-2))$ is not a subgraph of a graph contractible to $((0,0,0,0))$ and, thus,
it cannot serve as the dual graph of an SNC-curve with rational components contained in a smooth surface. 

\eexa


\section{Proof of the necessity part of the Main Theorem}\label{boundary}

\bnota\label{9.1} In this section we suppose that $\hX$ is an SNC-completion of a smooth rational semi-affine surface $X$.
As usual 
the dual graph of $\hD=\hX \setminus X$ will be denoted by $\Gamma$.
\enota

\bdefi\label{9.2} We say that an irreducible curve $F \subset X$ is distinguished if for any $\C$- or $\C^*$-fibration $f: X \to \C$
this curve is contained in a singular fiber of $f$.
\edefi

 Any automorphism of $X$ that transforms fibers of $f$ into fibers of $f$ respects general fibers (since only in neighborhoods of general fibers the fibration
is locally trivial). Hence,  the following fact holds.

\blem\label{lem6.5} Let $\nu$ be a complete vector field on $X$ that sends fibers of $f$ to fibers of $f$.   Then the flow of $\nu$ preserves every
singular fiber of $f$.
\elem

In combination with Theorem \ref{mthm1i} this leads to the following important technical tool which has already  appeared implicitly in Proposition \ref{nontoric}.

\bprop\label{9.3}
If a semi-affine surface $X$ contains a distinguished curve then its Remmert reduction $X_0$ is not generalized Gizatullin.
\eprop
 \bproof
Let $C$ be a distinguished curve and let $\nu$ be a complete algebraic vector field on $X_0$. Then by Theorem \ref{mthm1i} the flow of $\nu$ preserves a $\C$- or $\C^*$-fibration $f$. Since $C$ by definition belongs to a singular fiber of $f$ it is preserved by 
 Lemma \ref{lem6.5}. 
\eproof 

\bnota\label{9.4}
Recall that $\Gamma$ is contractible to a minimal graph $\Gamma_{\rm min}$ 
whose set of branch points $\mathrm{Br}(\Gamma_{\rm min}) $
is determined uniquely by Proposition \ref{reconstruction}.  Consider the subset $S$ of $\mathrm{Br}(\Gamma_{\rm min}) $
consisting of all branch points $E$ of valency 3 in $\Gamma_{\rm min}$ such that

\noindent $\bullet$ two of branches at $E$ are just $(-2)$-vertices;

\noindent $\bullet$ by elementary transformations (as in Proposition \ref{elementary}) on connected components of 
$\Gamma_{\rm min} \ominus \mathrm{Br}(\Gamma_{\rm min})$ one can make the weight of $E$ equal to $-1$. 

Note that $S$ is also determined uniquely by Proposition \ref{reconstruction} and Corollary \ref{cor.reconstruction}.
Let $T_{{\rm min}}=\mathrm{Br}(\Gamma_{\rm min}) \setminus S$ and  $\Gamma_{0,{\rm min}}, \Gamma_{1,{\rm min}}, \ldots , \Gamma_{n,{\rm min}}$
be connected components of $\Gamma_{\rm min} \ominus T_{{\rm min}}$. Let  $T$ be the set of the branch points in $\Gamma$ corresponding 
to $T_{{\rm min}}$ and $D_0,D_1, \ldots , D_n\subset \hD$ be the curves in associated with the connected components of $\Gamma\ominus T$.

\enota

Connectedness of $\hD$ and the description of vertices from $T_0$ in  Setting \ref{9.4}  imply the following.

\blem\label{9.5a} Each graph $\Gamma_{i, {\rm min}}$ has one of the following configurations:

{\rm (a)}  linear graph;

{\rm (b)} circular graph;

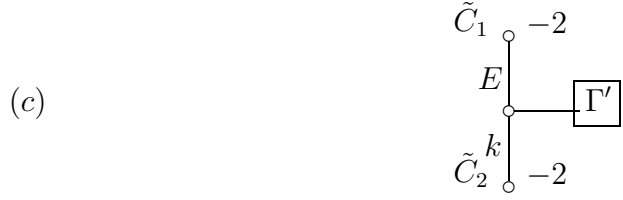
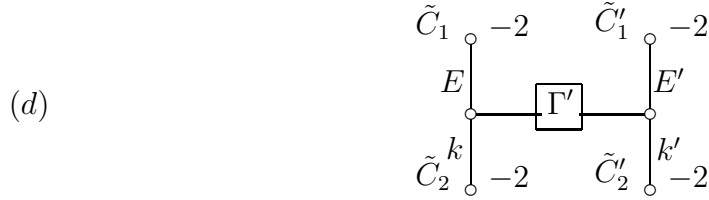
\vspace{2cm}

$$
\bigskip
\cou{E \, \, \, \, \, \, \,  }{k \, \, \, \, \, \, }\nlin\cshiftup{\tilde C_1}{-2}\slin\cshiftdown{\tilde C_2}{-2}\lin \Gamma' \mebox
\leqno{\,\,\,\,\,\,\, (c)}$$

\vspace{2cm}

$$
\bigskip
\cou{E \, \, \, \, \, \, \,  }{k \, \, \, \, \, \, }\nlin\cshiftup{\tilde C_1}{-2}\slin\cshiftdown{\tilde C_2}{-2}\lin \Gamma' \mebox\lin\cou{ \, \, \, \, \,
\, \, E'}{ \, \, \, \, \, \, \, k'}\nlin\cshiftup{\tilde C_1'}{-2}\slin\cshiftdown{\tilde C_2'}{-2}
\leqno{\,\,\,\,\,\,\, (d)}$$

\vspace{1cm}

\noindent where  $\Gamma'$ is a linear graph or empty. Furthermore,

\noindent $\bullet$ in (b) and (d) one has $\Gamma_{i,{\rm min}}=\Gamma_{\rm min}$;

\noindent $\bullet$ 
in (a) there are at most two neighbors of $\Gamma_{i,{\rm min}}$ in $\Gamma_{\rm min}$, each of them is contained in $T_{1,{\rm min}}$ and joined by an edge with an end-vertex of $\Gamma_{i,{\rm min}}$;

\noindent $\bullet$ in (c) there is at most one neighbor of $\Gamma_{i,{\rm min}}$ in $\Gamma_{\rm min}$, it is contained in $T_{1,{\rm min}}$ and joined by an edge with the right end-vertex of $\Gamma'$.
\elem

\blem\label{9.7} Let $\hat f : \hX \to \PP^1$  be a proper extension to $\hX$ of a $\C$- or $\C^*$-fibration $f: X \to \C$, i.e., some fiber  $F$
 of $\hat f$ has a support in $\hD$.
 Then this support is, in fact, contained in some $D_i$ (say, $D_0$). 

\elem
\bproof Making contractions in $\hD$ we get a pseudo-minimal extension 
$\hat g : \hX' \to \PP^1$ of $f$ with dual graph $\Gamma'$ for $\hD'=\hX' \setminus X$.
By Proposition \ref{fiber} $\Gamma'$ contains either  a linear 0-vertex $C$ or a subgraph of type $\Gamma_*$ which includes a chain $[[-2,-1,-2]]$ consisting of vertices $C_1+E+C_2$.

Since the image of $C$ in $\Gamma_{\rm min}$ is a linear vertex  it does not belong to $T_{{\rm min}}$,
i.e., the preimage of $C$ in $\hD$ (which is a support of a fiber of $\hat f$) is contained in $D_i$. 

In the case of the chain $C_1+E+C_2$ the images of $C_1$ and $C_2$ in $\Gamma_{\rm min}$ are linear vertices, 
while the image of $E$ is either a linear vertex of a nonnegative weight or a branch point of valency 3 which does not belong to $T_{{\rm min}}$.
Hence, the preimage of the chain in $\hD$ is contained in some $D_i$ which is the desired conclusion.

\eproof

\bprop\label{9.8}  Let $\hat f : \hX \to \PP^1$  and $\hat g : \hX \to \PP^1$ be proper extensions to $\hX$ of $\C$- or $\C^*$-fibrations on $X$.
Suppose that there is a fiber $F$ of $\hat f$ and a fiber $G$ of $\hat g$ with supports in
different $D_i$'s. Then $\hat f$ and $\hat g$ coincide (up to an automorphism of $\PP^1$) and, furthermore, every
$\PP^1$-fibration $\hat h : \hX \to \PP^1$ which extends a $\C$- or $\C^*$-fibration on $X$ coincides with $\hat f$.
\eprop

\bproof 
By assumption $F\cdot G =F^2=G^2= 0$ and for every section $S$ of $\hat f$ one has $F\cdot S=1$.
Hence, by the Hodge index theorem  $F$ and $G$ generate proportional elements $[F]$ and $[G]$ of $H^2(\hX, \Z )$, i.e., $n[F]=m[Q]$
for some nonzero $n,m\in \Z$. Note that $m=n$ since otherwise, say in the case of $m>n$, the product $G\cdot S$ is either negative or not an integer number.
Thus, these divisors $F$ and $G$ are equivalent
and $\PP^1$-fibrations $\hat f$ and $\hat g$ (generated by the linear systems of $F$ and $G$) coincide (up to an automorphism of the image $\PP^1$).
By Lemma \ref{9.7} $\hat h$ has a fiber with support in some $D_j$. Note that $D_j$ does not contain the support of either $F$ or $G$. Hence,
$\hat h$ coincides with either $\hat f$ or $\hat g$ and we are done.
\eproof

\bcor\label{9.9} Suppose that $\hX$ admits two distinct $\PP^1$-fibrations $\hat f: \hX \to \PP^1$ and $\hat g : \hX \to \PP^1$ that extend $\C$- or $\C^*$-fibrations on $X$.
Then both of them have fibers with support in the same $D_i$ (say, $D_0$) and no fibers with support in any $D_j$ where $j \ne i$.
\ecor

\bprop\label{9.10} If $X$ is generalized Gizatullin then distinct fibrations $\hat f$ and $\hat g$ as in Corollary \ref{9.9} exist.
In particular,   $\Gamma_{0,{\rm min}}$ cannot be a $0$-vertex.
\eprop

\bproof Assume that there is only one $\PP^1$-fibration $\hat f$ of this type. Then its restriction $f: X \to B:=f(X) \subset \C$ has no
singular fibers by Proposition \ref{9.3} and Lemma \ref{lem6.5}, i.e., $f$ is a locally trivial $\C$- or $\C^*$-fibration over $B$. Furthermore,
$B$ is not hyperbolic since otherwise for every complete vector field $\nu$ the
image of any integral curve under $f$ must be constant, i.e., $\nu$ must be tangent to the fibers of $f$ and $X$
has no open orbit. That is, $B$ is either $\C$ or $\C^*$. This implies that $X$ is isomorphic to either $\C^2$, $\C \times \C^*$, $(\C^*)^2$ or a twisted locally trivial $\C^*$-fibration over $\C^*$. In the first three cases our assumption is obviously wrong.

In the last case $\hD$ contains two fibers of $\hat f$ whose graphs under pseudo-minimality assumption are the chains $[[-2,-1,-2]]$ by Proposition \ref{fiber}. 
That is, $\Gamma$
is of the form (d) in Lemma \ref{9.5a} with $\Gamma'$  being the component $C_0$ of $\hD$ that is ramified double cover of $\PP^1$ under $\hat f$
and $k=k'=-1$. In order to avoid  the existence of another fibration   $\hat g$ one needs to require that $C_0^2 \leq -1$.
 However, in this case contracting the chains we get a Hirzebruch surface with the image of $C_0$ having self intersection at most 3 contrary
to Proposition \ref{dug.p1} which concludes the proof.
\eproof

\bexa\label{9.10a}  The last surface $S$ which is a twisted locally trivial $\C^*$-fibration over $\C^*$ is rather interesting.
Given a coordinate $z$ on a fiber $\C^*$ one can suppose that the monodromy around the puncture in the base in given by $z \to 1/z$.
Treating $\C^*$ as a complexification of a circle both in the base in the fiber one can see that $S$ is nothing but the complexification of the Klein bottle.
An SNC-completion of $S$ can be constructed in the following way. Consider the parabola $C$ given by $x-y^2=0$ in $\PP^1_x\times \PP^1_y$. 
Blow up point $(0,0)$ (resp.  $(\infty , \infty )$) and an infinitely near point. The resulting surface $\hat S$ is the desired completion of $S$
with fibers over $x=0$ and $x=\infty$ being $[[-2,-1,-2]]$-chains in $\hat S \setminus S$ and the proper transform $C_0$ of
$C$ playing the role of the  ramified double cover with $C_0^2=0$.  There are only two $\C^*$-fibrations $x$ and $y^2/x$ on $S$
corresponding to the $\PP^1$-fibrations associated with the chains $[[-2,-1, -2]]$ and vertex $C_0$ respectively.
Both $x$ and $y^2/x$ are rational first integrals (the first one for the complete algebraic field $\nu_1=(y^2-x ) \partial / \partial y$ and the second
one for $\nu_2 =2x \partial / \partial x + y\partial / \partial y$). There is also a complete algebraic field  $\nu =\nu_1+\nu_2$ \footnote{Completeness of $\nu$
can be extracted from completeness of $\nu_1$ and $\nu_2$ and the fact that $[\nu_1, \nu_2]=-\nu_1$.} for which neither $x$ nor $y^2/x$ is a rational first integral
(and which, therefore, has no rational first integral at all). Note that the function $x$ yields the preserved fibration on $S$. 

\eexa

\brem\label{D5}  Proposition \ref{9.10} does not hold in general for surfaces with open orbits.  Indeed, consider surfaces $yx^2+y^n+z^2=0, \, n\geq 4$ with
 du Val singularities of type $D_n$ mentioned in Example \ref{D4}. Each of them (say $S$) has a complete algebraic field tangent to the $\C^*$-fibration described in 
that example and the semi-simple field associated with the $\C^*$-action $(x,y,z) \to (\lambda^{n-1}x,\lambda^2 y,\lambda^{n}z)$. Hence, it has an open $\AAutH (S)$-orbit.
On the other hand there is an SNC-completion of $S$ with dual graph of the divisor at infinity as in  Proposition \ref{prop.2sur.rational} whose form in combination with 
Proposition \ref{fiber}
yields a unique $\C^*$-fibration and the absence of $\C$-fibrations on $S$. 

\erem

\bprop\label{9.11} In the Settings \ref{9.1} and \ref{9.4} let the Remmert reduction $X_0$ of $X$ be a generalized Gizatullin surface.
Then $T$ is empty and in particular $\Gamma_{\rm min}$ has one of the four configurations from Lemma \ref{9.5a}.
\eprop

\bproof By Lemma \ref{9.5a} it suffices to consider configurations (a) and (c) and show that $T$ is empty.

{\em Configuration (a).}   Assume that $\Gamma_{\rm min} \ne \Gamma_{0,{\rm min}}$, i.e., there is a vertex $V \in T_{{\rm min}}$ 
adjacent to the right end $E_1$ of $\Gamma_{0,{\rm min}}$.
By Proposition \ref{elementary} we can suppose that the left end $E_0$ of $\Gamma_{0,{\rm min}}$ is a $0$-vertex and by Proposition \ref{9.10}
$E_0 \ne E_1$. 

That is, the $\PP^1$-fibration $\hat f : \hX \to \PP^1$ induced by $E_0$ is an an extension

(i) of a $\C$-fibration if $E_0$ has no neighbor from $T_{\rm min}$; or

(ii) of a $\C^*$-fibration if $E_0$ has such a neighbor $W$ (where we allow equality $W=V$).

In (i) (and even in  (ii) when $W\ne V$) $V$ is not a section of the fibration $\hat f$ and is contained in a fiber $\hat f^{-1}(a)$ of $\hat f$. By
Lemma \ref{Cfeather} 
 $f^{-1} (a)$  is a singular fiber of $f$ that contains  $F=C \cap X \simeq \C$ where $C$ is a component of $\hat f^{-1} (a)$
that does not meet $D_0$.

Given a different $\C$- or $\C^*$-fibration $g: X \to \C$ by the Zariski theorem we can find an SNC-completion $\bX=X \cup \bD$ of $X$ which dominates $\hX$  and admits
a proper extension $\bar g$ of $g$. The preimage $\bar D_0$ of $D_0$ in $\bD$ has dual graph $\bar \Gamma_0$
which contains the dual graph $G_0$ of the preimage of $E_0$ and by Proposition \ref{9.8} the dual graph $G$ of the support of a fiber of $\bar g$. 
 Note that the proper transform $\bC$ of $C$
is again contained in the fiber of $\bar g$ because  it does not meet $\bD_0$.
By the Zariski lemma (e.g., see \cite[p. 90]{BPVdV}) $G_0 \setminus G$ cannot be empty.
This implies that either the proper transform of $V$ is not a section of the $\PP^1$-fibration $\bar g$ or  there are two neighbors of $G$ in $\bD$
which means that $g$ is a $\C^*$-fibration on $X$.  In both cases  $F \simeq \C$       is a component of a singular fiber of $g$ and, thus, it is a distinguished curve.
Since it survives the Remmert reduction 
we see that $X_0$ is not generalized Gizatullin by Proposition \ref{9.3}.  A contradiction, i.e., $T$ must be empty in this case.

When $W=V$ we deal only with $\C^*$-fibrations. Blowing up the edge between $E_1$ and $V$, if necessary,
 replace $\Gamma_{0,{\rm min}}$
with a linear graph $\Gamma_0$ which contains at least three vertices.  Using reconstructions
as in  Proposition \ref{elementary} we can change $\Gamma_0$ so that it contains now a $0$-vertex $E$ somewhere in the middle (i.e., $E$ is not a neighbor of $V$).
Lemma \ref{Cfeather} supplies us as before with a distinguished curve $F \simeq \C$ which concludes the consideration of configuration (a).

{\em Configuration (c)}.  Consider  $\hat f : \hX \to \PP^1$ induced by the chain $\tilde C_1+E+\tilde C_2=[[-2,-1,-2]]$ contained in $D_0$.
  Note that $D_0$ contains an irreducible component $U$ meeting $E$ transversely. 
Since the multiplicity of $E$ in the fiber of $\hat f$ is two the map $\hat f|_U : U \to \PP^1$ has degree two.
This implies  that $ f$ is a $\C^*$-fibration on $X$. If there exists $V\in T$ then $V\ne U$ since otherwise the second fibration $\hat g$ does not exists. 
Hence, by
Lemma \ref{Cfeather}  we have again  a singular fiber of $f$ containing an irreducible component $F$ isomorphic to $\C$
and the same argument shows that $g$ as before is also a $\C^*$-fibration with $F$ being a component of a fiber of $g$. This yields
a distinguished curve and the desired conclusion.

\eproof

Now we are able to prove one direction of the Theorem \ref{mthm}, namely the necessary condition for a semi-affine surface to be generalized Gizatullin.

\bprop\label{prop7.1} Let $X$ be a semi-affine surface which is generalized Gizatullin, then $X$ possesses a
completion $\bar X$ where the dual graph of $\bar X\setminus X$
is of one of the forms (1) - (6) from Theorem  \ref{mthm}.
\eprop

\bproof  We go through the cases of Lemma \ref{9.5a}.

 For the case of a linear dual graph the Proposition is well known.  Let us consider first circular graphs.  By \cite[Proposition 3.28]{FKZ} they can be reduced via birational transformations to
the following standard (and ``essentially unique") forms:

(i) $((0_{2k}, w_1, \ldots , w_n))$ with $k\geq 0, n>0$, and $w_i \leq -2$; or 

(ii) $((0_l, w))$ with $l >0$ and $w\leq 0$; or 

(iii) $(0_{2k}, -1,-1))$ with $k\geq 0$

\noindent where the subindex reflects the number of consequent zero weights. Note that in (i) we disregard the case of $k=0$ (since the intersection matrix of such
a graph is negative definite contrary to the Nakai-Moishezon criterion) and $k\geq 2$ (to avoid a contradiction with the Hodge index theorem).
Similarly, in (ii) we omit the case of $l=1$ and $w \leq -1$ because otherwise we get only one $\PP^1$-fibration as in Proposition \ref{9.9} contrary to Proposition \ref{9.10}.
The Hodge index theorem implies also that in (ii) $l \leq 3$. Hence, the remaining possibilities in (i) and (ii) produce (2a) and (2b). By the same arguments in (ii) we have to
consider only weights as in (2c). 

By Proposition \ref{9.11}
it suffices now to consider $\Gamma_{\rm min}$ as in configurations (c) and (d) in Lemma \ref{9.5a} and show that $\Gamma'$ is of the desired form.

Consider first the case when there exists  a non negative vertex in a minimal graph of $\Gamma'$. After blowing up which
keeps the graph linear  we can always suppose that it is actually of weight $0$. Using operations of  of form $[[v,0,w]] \to [[v-1,0, w+1]]$ as in 
Proposition \ref{elementary} we can suppose furthermore, that this vertex is the left endpoint $C_0$ of $\Gamma'$ and the weight $k$ in configurations (c) and (d) from Lemma \ref{9.5a}  is $-1$.   Proposition \ref{dug.p0} implies now that $\Gamma$
must be of desired form which describes configurations  (3) and (5)  in Theorem \ref{mthm} completely.

Now consider the case when there is no non negative vertex in a minimal graph of $\Gamma'$. If in (d) $\Gamma'$ is not empty then
the divisor $$C_1 + \ldots + C_n + \tilde C_1' + E' + \tilde C_2'$$ is contained in the same fiber
of $\hat f$.  Hence,  by the Zariski lemma (e.g., see \cite[p. 90]{BPVdV}) we have not only $C_i^2\leq -2$
but also  $(E')^2 \leq -1$ since the intersection matrix of this divisor must be negative definite. 
By Corollary \ref{cor.reconstruction} 
we see that  this minimal graph is unique and $\hat f$ is the only extension of a $\C$- or $\C^*$-fibration on $X$ contrary to Proposition \ref{9.10}.
Thus, $\Gamma'$ is empty in this case and we have (6) (condition $k' \geq -1$ is necessary to provide a second $\PP^1$-fibration not equal to $\hat f$). 

In the absence of a non negative vertex for (c) we have $k=-1$ and unless $w_0=-2$ there is again the same contradiction with  Proposition \ref{9.10}.
Thus, $w_0=-2$ and we have three $\PP^1$-fibrations associated with with the chains $\tilde C_1 + E + \tilde C_2$, $\tilde C_1 + E + C_0$, and $\tilde C_2 + E + C_0$
respectively.
However, if $n \geq 1$ the restriction of any of the last two fibrations to $X$ is not a $\C^*$-fibration because contrary to Proposition \ref{fiber}(3b) the $\tilde C_i + E + C_0$
meets not only the ramified double cover but also $C_1$. This leads to case (4) and we are done.
 \eproof

\section{Proof of the sufficiency part of the Main Theorem}\label{proof}

\blem\label{lem8.2} Let $X$ be a smooth semi-affine surface. Suppose that $f_i : X \to \C , \, i=1,2$ are either $\C$ or $\C^*$-fibrations  such that the intersection of every pair of non-singular fibers of $f_1$
and $f_2$ is a finite non-empty set.
Let $S_i$ be the union of singular fibers of $f_i$. Then there is an open orbit $U$ of the natural $\AAutH (X)$-action (Definition \ref{int.2} )
on $X$ such that its complement is contained in $S_1 \cap S_2$.
\elem

\bproof Let $x \notin S_1\cap S_2$ and let $U$ be the orbit of $x$. Say $x$ is contained in a non-singular fiber $E_1$ of $f_1$. By Lemma \ref{tangent.field}
$E_1\subset U$,
i.e.,  we can suppose that $x$ is an arbitrary point of $E_1$.  In particular, we can suppose that $x \in E_1 \cap E_2$ where $E_2$ is
a given non-singular fiber of $f_2$.
By Lemma \ref{tangent.field} $E_2\subset U$. Similarly any given fiber of $f_1$ is contained in $U$ and we are done.
\eproof

\proof[Proof of Theorem \ref{mthm}] By Proposition \ref{prop7.1} it suffices to show that every normal affine algebraic surface $X_0$ with a dual graph of the boundary appearing in Theorem \ref{mthm} is a generalized Gizatullin surface.

If $X_0$ is a Gizatullin surface then it is known that $X_0$ is already quasi-homogeneous under the algebraic automorphisms by \cite{Gi}. Thus, we need to show
now that cases (2)-(6) in Theorem \ref{mthm}  indeed present surfaces quasi-homogeneous under the group of  algebraically generated automorphisms. Suppose that $X\to X_0$ is
a minimal resolution of singularities, i.e., $X$ is smooth semi-affine. Consider 
the most difficult dual graph $\Gamma$ of $\bD =\bX \setminus X$ as in Figure (5). Making inner blowing up if necessary we can suppose that the weight of $C_0$ is zero.
Let $f: X \to \C$ be the twisted $\C^*$-fibration associated with the subgraph $K=\tilde C_1 +E + \tilde C_2$ and $f_0: X \to \C$ be the
untwisted $\C^*$-fibration associated with the 0-vertex $C_0$.
Suppose that $\bar f : \bX \to \PP^1$ and $\bar f_0 : \bX \to \PP^1$ are their proper extensions. Let $S$ (resp $S_0$)
be the union of singular fibers of $\bar f$ (resp. $\bar f_0$) that meet $C_0$ (resp. $ K$) and $S'$ (resp. $S_0'$) be the
union of such fibers that do not. That is, each fiber from $S'$ (resp. $S_0'$) meets one of curves
$C_1, \ldots , C_n, E', \tilde C_1', \tilde C_2'$.  Let $U$ be the open orbit of the natural $\AAutH (X)$-action in $X$. By Lemma \ref{lem8.2}
$X \setminus U$ is contained in $(S\cup S') \cap (S_0 \cup S_0')$. On the other hand, since $S$ is not contained in a fiber of $\bar f_0$,
with an exception of a finite set $S$ is contained in $U$ by Lemma \ref{tangent.field}.  The same is true for $S_0$. Thus, $X \setminus U$ is contained in $(S' \cap
S_0') \cup T_0$
where $T_0$ is a finite set. That is, we need to show that up to a finite set  every curve $F\subset S' \cup S_0'$ that is a component of singular fibers of both $\bar f$ and
$\bar f_0$ is contained in $U$.
Let $F$ meet $C_1$. Then by Proposition \ref{elementary} after a sequence of elementary transformations such that all vertices of $\Gamma$
but $C_0$ survive we can make the weight of $C_1$ equal to $0$. In particular, this $0$-vertex yields a $\C^*$-fibration $f_1: X \to \C$.
Since $F$ meets $C_1$ we see that $F\cap X$ is not contained in a fiber of $f_1$ and, therefore, by Lemma \ref{tangent.field} one has $F \setminus T_1 \subset U$
where $T_1$ is a finite set. 

Similarly, by Proposition \ref{elementary}  if $i \geq 2$ then after a sequence of elementary transformations under which all vertices 
$C_i, \ldots , C_n, E', \tilde C_1', \tilde C_2'$ survive we can make the weight of $C_i$ equal to $0$. Thus, the same argument
implies that if $F$ meets $C_i$ then  $F \setminus T_i \subset U$ where $T_i$ is a finite set. In the case of $F$ meeting of the curves
$E', \tilde C_1'$, or  $\tilde C_2'$ using elementary transformations in the chain $C_0+ \ldots +C_n$ we can make the weight of $E'$ equal to
$-1$. Then $\tilde C_1' + E' + \tilde C_2'$ becomes a subgraph that induces a twisted $\C^*$-fibration whose restriction to
$F\cap X$ is not constant. That is, up to a finite set $F$ is contained in $U$. Hence, we are done in case (5).

The argument for case (3) and a circular graph in case (2) are similar (say, the only difference in (3) is that when one makes the weight of $C_n$ equal to 0
then the associated fibration is a $\C$-fibration and not a $\C^*$-fibration).

Also in case (4) the argument is similar, we have to work with the three twisted $\C^*$-fibrations associated with the three $[[-2,-1,-2]]$ subgraphs.

If $k'\geq 0$ in case (6) then one needs to make a sequence of inner blow-ups over the edge between $E$ and $E'$ such that the resulting graph looks like\\

$$\bigskip
\Gamma=\hspace{1cm}
\cou{E \, \, \, \, \, \, \,  }{-2 \, \, \, \, \, \, \, \, \, }\nlin\cshiftup{\tilde C_1}{-2}\slin\cshiftdown{\tilde
C_2}{-2}\lin\cou{C_0}{-2}\lin\cou{C_1}{-2}\lin\ldots
\lin\cou{C_{n-1}}{-2}\lin\cou{C_{n}}{-1}\lin\cou{ \, \, \, \, \, \, \, E'}{ \, \, \, \, \, \,\, \,\,-1}\nlin\cshiftup{\tilde C_1'}{-2}\slin\cshiftdown{\tilde
C_2'}{-2}
\quad .$$\\
Then we have two twisted $\C^*$-fibrations $g': X \to \C$ and $g : X \to \C$ induced by the subgraphs $K' =\tilde C_1' + E' + \tilde C_2'$ and $K=\Gamma\ominus
K'$, indeed $K$ is contractible to a $[[-2,-1,-2]]$ subgraph by contracting $C_0+\ldots + C_n$. Suppose that $\bX$ is the SNC-completion of $X$ with
the boundary described by the graph above and $\bar g : \bX \to \PP^1$ (resp.  $\bar g' : \bX \to \PP^1$) is a proper extension of $g$ (resp. $g'$).
Note that the singular fibers of $\bar g$ must meet $ K'$ but not $ K$ while for the singular fibers of $\bar g'$ the situation is reversed.
In particular, only complete curves that are contained in  $X$ may be common components of singular fibers of $\bar g$ and $\bar g'$. 
By Lemma \ref{lem8.2} $U$ is contained in the complement to the union of such components in $X$.
Hence, $X_0$ is a generalized Gizatullin surface since these components are contractible to points in the Remmert reduction.

In the case of $k'=-1$ a similar argument works and we are done.

\section{Homogeneity}\label{homogen}

\bnota
In this section $X$ is a smooth affine surface with an SNC-completion $\bX$ such that the dual graph of $\bD = \bX \setminus X$ is one of those in Theorem \ref{mthm}.
In particular, $X$ is a generalized Gizatullin surface.
\enota

Note that if $X$ admits a $\C$ or $\C^*$-fibration without singular fibers (say, $X$ is the complexification of the Klein bottle from Example \ref{9.10a}) then it is homogeneous with respect to the  $\AAutH(X)$-action
because of the absence of fixed points for this action by virtue of Lemma \ref{tangent.field}.
The same remains true in several other cases.

\bthm\label{circular}
Let $\bD$ have a circular dual graph as in (2) of Theorem \ref{mthm}. Then $X$ is homogeneous with respect to the natural  $\AAutH(X)$-action.
\ethm

\bproof
Suppose that  $C_0, \ldots , C_n$ are irreducible components of $\bD$ such that $C_0^2=C_1^2 = 0$ and $C_i^1\leq -1$ for $i\geq 2$,
$C_i\cdot C_j= 1$ for $\vert i-j\vert = 1$ or $\lbrace i,j \rbrace = \lbrace 0,n \rbrace$, and $C_i\cdot C_j = 0 $ otherwise.
Let $\bar f : \bX \to \PP^1$ be the $\PP^1$-fibration associated with $C_0$, the $\C^*$-fibration $f$ be its restriction to $X$ and let
\{$F_i \}$ be the irreducible components of singular fibers of  $f$. 

{\em Case 1}: $n\geq 3$.  The absence of branch points in the dual graph of $\bD$ and
the smoothness of $X$ imply that all singular fibers of $\bar f$ but the one 
containing $C_2$ (say $\bar f^{-1}(0)$) are
chains $[[-1,-1]]$ while $\bar f^{-1}(0)$  consists of the chain $\cC = C_2 + \ldots + C_{n-1}$ (joining sections $C_1$ and $C_n$) 
and some other components
adjacent to smooth points of this chain which are $(-1)$-curves because of Lemma \ref{irreducible}. Each of these $(-1)$-curves is of course a
closure $\bar F_i$ of some $F_i$ and it is a component of a singular fiber of $\PP^1$-fibration associated with $C_1$.

Hence, by Lemma \ref{lem8.2} one can suppose that a potential fixed point of the $\AAutH (X)$-action is contained not only in a chain $[[-1,-1]]$ mentioned before but
also in some $F_i$ from $f^{-1}(0)$. Therefore,
it is enough to show that for any given $F_i$ there is a complete algebraic vector field whose restriction to $F_i$ is locally nilpotent and nontrivial, i.e.,
it generates a translation on $F_i$. Though a priori $F_i$ may be adjacent to any $C_j$ with $2\leq j \leq n$
a reconstruction as in Proposition \ref{elementary} enables us to consider only the case when $F_i$ is adjacent to $C_2$. 

Contracting in fibers of $\bar  f$ irreducible component not adjacent to $C_1$ (in particular $C_2$ is not contracted) we get a morphism $\varphi : \bX \to \bX'$ into a 
Hirzebruch surface  $\bX'$ with $C_1'$ and $C_n'$ playing the roles of disjoint sections where $C_i'$ is the proper transform of $C_i$ in $\bX'$. That is, $\bX' \setminus C_0'$
is naturally isomorphic to $\C_x \times \PP_y^1$ with  $C_1'\setminus C_0',C_2'\setminus C_0',C_n'$ given by $\lbrace y=\infty \}$,
$\lbrace x=0 \}$, $\lbrace y=0 \}$. 

Lemma \ref{field} implies now that the  pull-back of  the vector field $\mu=y\frac{\partial}{\partial y}$ on $\bX'\setminus C_0'$ is  a rational vector field $\bar \mu$ on $\bar X\setminus C_0$ which has 
only simple poles and they are located on those $\bar F_i$'s that are adjacent to the chain $\cC$. This means that $x\mu$  induces  a regular vector field $\nu$ on $X$
and even on $\bX \setminus C_0$.

Note that $\varphi (F_i)=(0,y_0) \in  \C^* \times \PP^1$ with $y_0 \ne 0, \infty$ and $F_i$ is obtained as the result of a monoidal transformation at this point.
That is, one can introduce a local coordinate system $(u,v)$ on $X$ such that $\varphi (u,v)=(u,uv+y_0)$ and $F$ is given by equation $u=0$.
Then $\nu$ is given by $(uv+y_0)\frac{\partial}{\partial v}$, i.e., its restriction to $F_i$ is nonzero and locally nilpotent.
Hence, no point of $F_i$ is fixed under the $\AAutH(X)$-action which implies the desired conclusion in this case.

{\em Case 2}: $n=2$.  One can blow up the edge between $C_1$ and $C_2$ to get an extra vertex $C$, i.e., we have four vertices in the new dual graph. 
Consider $\bar f$, $f$, $F_i$, and $\bar f^{-1}(0)$ as before. Note that
the weight of the proper transform of $C_1$ becomes $-1$ but elementary transformation from Proposition \ref{elementary} can make
it again zero while keeping $C$ intact. That is, any $F_i$ from $\bar f^{-1}(0)$ is contained in a singular fiber of a $\C^*$-fibration on $X$ different from $f$.
Lemma \ref{lem8.2} implies that it suffices again to construct a translation on $F_i$ and the previous argument works. 

{\em Case 3}: $n=1$ and $\bD$ consisting of two $0$-components $C_0$ and $C_1$ meeting each other transversely at two points. It requires a different approach
which we sketch below. Let $\bar f$ be again the $\PP^1$-fibation on $\bX$ associated with $C_0$. Making contraction $\varphi : \bX \to \bX'$ in the fibers of $\bar f$
we get a Hirzebruch surface $\bX'$ with the proper transform $C_1'$ of $C_1$ playing the role of a ramified double section. 
Let $\bar g : \bX' \to \PP^1_x$ be induced by $\bar f$.  Without loss of generality
we can suppose that $C_0'=\bar g^{-1}(1)$ while $0$ and $\infty$ are the only singular values of $\bar g |_{C_1'}$ as in Remark \ref{lem.2singfib}. Furthermore,
applying the same reconstruction we exploited in Proposition \ref{9.10} we can suppose that 
$C_1'$ is the closure of the curve given by $x=y^2$ in $\bX' \setminus \bar g^{-1} (\infty ) \simeq \C_x\times \PP^1_y$.
That is, the restriction $g$ of $\bar g$ to $X'=\bX'\setminus (C_0'\cup C_1')$ is a $\C^*$-fibration with two singular fibers $g^{-1} (0)$ and $g^{-1}(\infty )$ (both isomorphic to $\C$).

Since $X$ does not contain complete curves the surface $\bX$ is obtained from $\bX'$ by several monoidal transformations at different points of $C_1'$. 
Hence, the singular fibers of $f$ are $f^{-1}(0)$, $f^{-1}(\infty )$ (each of them consists of one
or two connected components isomorphic to $\C$) and fibers that are unions of form $F_1\cup F_2$ where $F_i \simeq \C$ and $F_1$ meets $F_2$ transversely at one point.

The vector field $ \frac{(x-y^2)}{x-1} \frac{ \partial}{ \partial y}$
is regular and complete on $X'$ and its restriction to $g^{-1} (0)$ and $g^{-1}(\infty )$ 
induces nontrivial translations.  Furthermore, a calculation shows that it induces a regular vector field $\nu$ on $X$ which is a translation on every irreducible component of $f^{-1}(0)$ or $f^{-1}(\infty )$.
Therefore, by Lemma \ref{field} points of type $F_1 \cap F_2$ are the only potential  fixed points of the $\AAutH (X)$-action. 

Consider the following reconstruction of the boundary divisor: Blow up one of edges between $C_0$ and $C_1$ and contract the proper transform of $C_1$.
The resulting completion $\hX$ of $X$ has the boundary divisor $\hD$ consisting of two $0$-vertices $\hat C_0$ and $\hat C_1$ where $\hat C_0$ is the proper
transform of $C_0$. Let $\hat f : \hX \to \PP^1$ be the $\PP^1$-fibration on $\bX$ associated with $\hat C_0$ such that $\hat C_0=\hat f^{-1}(1)$.
Consider the fibers $\hat F_1 \cup \hat F_2$ of $\hat f|_X$ similar to $F_1 \cup F_2$, i.e., every point that is not of type $\hat F_1 \cap \hat F_2$
belongs to the open orbit $U$ of the $\AAutH (X)$-action. Note that by construction the proper transform $G_i$ of $F_i$ meets $\hat C_0$ transversely at one point.
Hence, its intersection with every other fiber of $\hat f$ is also 1. In particular, $G_i$ cannot meet the fiber $\hat F_1 \cup \hat F_2$ at the double point $\hat F_1 \cap \hat F_2$.
This implies that $G_1 \cap G_2 \ne \hat F_1 \cap \hat F_2$. Thus, $F_1 \cap F_2 \in U$ and we are done.
\eproof

Recall that there are Gizatullin surfaces that are not homogeneous with respect to the natural $\Aut$-action. A list of such surfaces appeared in \cite{Ko} and
we show that every surface from  this list is homogeneous with respect to  the natural $\AAutH$-action. In fact this is true for a wider collection of
Gizatullin surfaces to describe which we need to remind the following. 

Let $\bY$ be an SNC-completion of a smooth Gizatullin surface $Y$ by a standard zigzag $\bD =\bY \setminus Y  =C_0 \cup \ldots \cup C_{n-1}, \, n\geq 3$. 
The $0$-vertices $C_0$ and $C_1$ of the zigzag
induce two $\PP^1$-fibrations that
lead to a morphism $\bar \Phi =(\bar \varphi_0 , \bar \varphi_1) : \bY \to \PP^1\times \PP^1$ with restriction
$\Phi =(\varphi_0, \varphi_1) : Y \to \C_{x,y}^2$. Omitting a simple case of $n=3$ we suppose further that $\bar \Phi (C_3 \cup \ldots \cup C_{n-1})=(0,0)$,
i.e.,  the only singular fiber $\bar \varphi_0^{-1} (0)$ of $\bar \varphi_0$ is contracted by $\bar \Phi$ to the proper transform of $C_2$.  
The  components of  $\bar \varphi_0^{-1} (0)$ different from $C_2 \ldots , C_{n-1}$ are called feathers (in terminology of \cite{FKZ.a} or \cite{Ko}). 
For every surface in Kovalenko's list each feather is a $(-1)$-curve.

\bthm
Let $Y$ be a smooth Gizatullin surface $Y$ such that every feather is  a $(-1)$-curve. Then $Y$ is homogeneous with respect to  the natural $\AAutH$-action.
\ethm

\bproof  
Since each feather is a $(-1)$-curve they can be contracted first.
This implies that for the sequence $\bar \Phi : \bY \to \PP^1\times \PP^1$  of monoidal transformations, $C_{n-1}$ is obtained from the proper transform $0 \times \PP^1$ of $C_2$
after several (say $k$) outer
blowing-ups in   (see Section \ref{dualgraphs} for definition of outer blowing up) at the origin and infinitely near points. 
Hence, for some fixed values $a_1, \ldots , a_{k-1} \in \C$ and general $b \in \C$ the proper transform $C$ of the curve $y=a_1x+ \ldots + a_{k-1}x^{k-1} +bx^k=: h(x)$
in $\bY$ meets $C_{n-1}$ at a general point. The triangular automorphism $(x,y) \to (x, y-h(x))$ of $\C^2$ induces an isomorphism of $Y$ on another Gizatullin surface
$Y'$ which has a completion $\bY'$ by a standard zigzag $C_0' + \ldots + C_{n-1}'$ such that this isomorphism extends regularly to $\bY \setminus C_0 \to \bY' \setminus C_0'$.
We replace $\bY$ by $\bY'$. The advantage is that $C$ is now the proper transform of the $x$-axis in $\C^2$, i.e., it meets
both $C_{n-1}$ and $C_0$.
That is, the graph of $\bD \cup C$ becomes  circular with $C$ playing the role of $C_n$. Thus,
$X=Y \setminus C$ is a surface of type (2) from Theorem \ref{mthm} and by Theorem \ref{circular}
$\AAutH(X)$ acts transitively on $X$.  

Recall that the field $\nu$ in Case 1 of Theorem \ref{circular} extends regularly to $\bX \setminus C_0$ and in particular to $C_n\setminus C_0$ (and, therefore, to $Y$).
Furthermore, consider transformations $((0,0,w_3, \ldots ,w_n))\to ((w_3, \ldots ,w_{j-2},0,0, w_j, \ldots , w_n)) $
from Proposition \ref{elementary} used in Case 1 to make a feather adjacent to $C_2$ instead of $C_j$.  
Note that $C_n$ survives such a transformation and plays the role of $C_{n-j+3}$ in the modified graph,
i.e., it is still contained in $\bX \setminus C_0$. 
Hence,  even after these transformations the  flow of $\nu$ is
extendable to $Y$. Since the homogeneity of $X$ is provided by elements of these  flows we see that $X$ is contained actually in the open orbit of the $\AAutH(Y)$-action.
Note that  $C\cap Y$ does not contain fixed points of the $\AAutH(Y)$-action since each point of $C\cap Y$
can be moved by a $\G_a$-action  induced on $Y$
by the field of form $x^m \frac{\partial}{\partial y}$ on $\C^2$ with $m >>0$.
Thus, the open orbit coincides with $Y$. Therefore, $Y$ is, indeed,  $\AAutH(Y)$-homogeneous. 
\eproof

\providecommand{\bysame}{\leavevmode\hboxto3em{\hrulefill}\thinspace}

\end{document}